\documentclass[journal,draftcls,onecolumn,letterpaper,twoside,11pt]{IEEEtran}

\usepackage[margin=1in]{geometry}
\usepackage{amssymb, mathrsfs,verbatim} 
\usepackage[cmex10]{amsmath}
\usepackage{cite}
\RequirePackage[colorlinks,linkcolor=black,citecolor=black,urlcolor=blue]{hyperref}
\usepackage{enumerate}
\usepackage[dvips]{graphicx}
\usepackage{verbatim}
\usepackage{bbm}
\usepackage{subfig}
\usepackage{epsfig}
\usepackage{epstopdf}
\usepackage{graphicx} 

\usepackage{float,latexsym,amsfonts,amstext,times}

\usepackage{chngcntr}
\usepackage{apptools}

\usepackage{color}

\newtheorem{theorem}{Theorem}[section]

\newtheorem{proposition}{Proposition}[section]
\newtheorem{corollary}{Corollary}[section]
\newtheorem{lemma}{Lemma}[section]

\newcommand{\cS}{\mathcal{S}}
\newcommand{\cJ}{\mathcal{J}}
\newcommand{\cD}{\mathcal{D}}
\newcommand{\cC}{\mathcal{C}}
\newcommand{\cF}{\mathcal{F}}
\newcommand{\cP}{\mathcal{P}}

\newcommand{\Pro}{\mathsf{P}}
\newcommand{\Exp}{\mathsf{E}}
\newcommand{\bN}{\mathbb{N}}

\newcommand{\bS}{\mathbb{S}}

\begin{document}

\title{ Sequential anomaly detection  \\under sampling constraints}
\author{Aristomenis Tsopelakos and  Georgios Fellouris, 
	\textit{Member, IEEE}}

\thanks{Georgios Fellouris is with the Department of Statistics, the Electrical and Computer Engineering Department, and the Coordinated Science Laboratory,
	University of Illinois at Urbana–Champaign, Champaign, IL 61820 USA, (e-mail: fellouri@illinois.edu). Aristomenis Tsopelakos is with  the Electrical and Computer Engineering Department and the Coordinated Science Laboratory, University of Illinois at Urbana-Champaign, Urbana IL 61801, USA, (e-mail: tsopela2@illinois.edu). 
	The work of the two authors  were supported in part by the NSF under Grant CIF 1514245, through the University of Illinois at Urbana–Champaign. The work of  Georgios Fellouris was also supported in part by the NSF under Grant DMS 1737962, through the University of Illinois at Urbana–Champaign.}
\maketitle

\begin{abstract}
The problem of sequential anomaly detection is considered, where multiple data sources are monitored in real time  and  the goal is to identify the ``anomalous'' ones among them,  when it is not possible to sample all sources at all times.  A  detection scheme  in this context requires specifying not only   when to stop sampling and which sources to identify as anomalous upon stopping, but also which sources to sample at each time instance until stopping.   A novel formulation  for this problem  is proposed, in which the  number of anomalous sources is not necessarily known in advance and the number of sampled sources per  time instance is not necessarily fixed. Instead, an arbitrary lower bound and an arbitrary upper bound  are assumed on the number of anomalous sources, and the fraction of the  expected number of samples over the expected time  until stopping is required to  not exceed  an arbitrary, user-specified level. In addition to this sampling constraint,  the probabilities of at least one false alarm and at least one missed detection are controlled below  user-specified tolerance levels.  A general criterion  is established  for a policy to achieve the minimum expected time until stopping to a first-order asymptotic approximation  as  the two familywise error rates go to zero. Moreover, the asymptotic optimality is established of a family of policies that sample each source at each time instance  with a  probability that depends on  past observations only through the current estimate of the subset of anomalous sources. This family includes, in particular, a novel policy that requires minimal computation  under any setup of the problem.
\end{abstract}


\section{Introduction}\label{sec:Intro} 
In various engineering and scientific areas data are often  collected   in real time over  multiple streams, and  it is of interest to quickly  identify  those data streams, if any,  that exhibit outlying  behavior.  In brain science, 
for example,  it is desirable to  identify  groups of cells with  large  vibration frequency, as this is a symptom for the development of a particular malfunctioning \cite{Brain_2010}. In fraud prevention security systems in e-commerce, it is desirable to identify transition links with low transition rate, as this may be an indication that a link is tapped  \cite{Fraud_2002}.
Such applications, among many others,  motivate the study of sequential multiple testing problems where  the data for the various  hypotheses are generated  by distinct sources, there are two hypotheses for each data source, and the goal is to identify as quickly as possible  the ``anomalous'' sources, i.e., those   in which the alternative hypothesis is correct.   In certain works, e.g., \cite{De_and_Baron_2012a,De_and_Baron_2012b,Bartroff_and_Song_2014,Song_and_Fellouris_2016,Song_and_Fellouris_2019,Bartroff_and_Song_2020}, 
it is  assumed that all sources are sampled at each time instance, whereas in others, e.g.,    \cite{Zigangirov_1966, Klimko1975OptimalSS,Dragalin_1996, Cohen2015active,huang2017active,oddball_2018,Cohen2019nonlinearcost,Tsopela_2019,Cohen2020composite},  
only a fixed  number of sources  (typically, only one)  can be sampled at  each time  instance.  In the latter case, apart from  when to stop sampling  and  which data sources  to  identify as anomalous upon stopping, one also needs to specify  which  sources to sample at every time  instance until stopping.

The latter problem, which is often called  \textit{sequential anomaly detection} in the literature,  can  be viewed as a special case of the  \textit{sequential  multi-hypothesis testing problem with  controlled sensing (or observation control)}, where  the  goal is to solve a sequential multi-hypothesis testing problem  while taking  at each time  instance an action that influences the distribution of the  observations \cite{chernoff1959,albert1961,Bessler1960_I,Bessler1960_II,Kiefer_Sacks_1963,Keener_1984,Lalley_Lorden_1986,nitinawarat_controlled_2013,naghshvar2013active, nitinawarat_controlled_2015,Aditya_2021}. 
In the anomaly detection case, the action is the selection of the sources to be sampled, whereas  the  hypotheses correspond to the possible subsets of anomalous sources. Therefore, policies and results in the context of sequential multi-hypothesis testing with controlled sensing  are applicable, in principle at least,  to the  sequential anomaly detection problem.  Such a policy was first proposed in   \cite{chernoff1959} in the case of two hypotheses,  and subsequently generalized in  \cite{Bessler1960_I,Bessler1960_II} to the case of an  arbitrary,  finite number of hypotheses.   When applied to the sequential anomaly detection problem, this policy samples  each subset of sources of the allowed size at each time instance with a certain probability that depends on the past observations only through the currently estimated subset of anomalous sources. 

In general, the implementation of the policy in  \cite{Bessler1960_I} requires solving,  for each subset of anomalous sources, a linear system where the number of equations is equal to the number of sources and the number of unknowns is equal to the number of all subsets of sources of the allowed size. Moreover, its asymptotic optimality  has been established only under  restrictive assumptions, such as  when the following hold simultaneously: it is known \textit{a priori} that there is only \textit{one} anomalous source,  it is possible to  sample only \textit{one} source at a time,  the  testing problems are identical, and the sources generate  \textit{Bernoulli} random variables under each hypothesis \cite[Appendix A]{Bessler1960_II}.  To   avoid  such restrictions, it has been proposed to modify the policy in \cite{Bessler1960_I} at an appropriate  subsequence of time instances, at which  each  subset of sources of the allowed size is sampled with the same probability \cite[Remark 7]{chernoff1959},\cite{nitinawarat_controlled_2013}. Such a \textit{modified} policy  was shown in  \cite{nitinawarat_controlled_2013} to always  be asymptotically optimal, as long as the log-likelihood ratio statistic of  each observation has a  finite \textit{second} moment.

A goal of the present work is to show that  the \textit{unmodified} policy in \cite{Bessler1960_I}   
is always asymptotically optimal  in the context of the above sequential anomaly detection problem, as long as the log-likelihood ratio statistic of each observation has a  finite \textit{first} moment. However, our main goal in this paper is  to propose a more general framework for the problem of sequential anomaly detection with sampling constraints  that  (i) does not rely on the restrictive assumption that the number of anomalous sources is known in advance,  (ii) allows for two distinct error constraints and captures the asymmetry between a \textit{false alarm}, i.e.,  falsely identifying a source as anomalous,  and a \textit{missed detection}, i.e., failing to detect an anomalous source, and most importantly,  (iii) relaxes the hard sampling constraint that  the same  number of sources must  be sampled  at each time instance, and  (iv) admits an asymptotically optimal solution that is  convenient to   implement  under any setup of the problem.

To be more specific, in this paper we  assume  an arbitrary, user-specified lower bound and an arbitrary, user-specified upper bound on the number of anomalous sources.  This setup includes  the case of no prior information, the case where  the number of anomalous sources  is known in advance, as well as more realistic cases of prior information, such as when there is  only a non-trivial upper bound on the number of anomalous sources.   Moreover, we require control of the probabilities of at least one false alarm and at least one missed detection below arbitrary, user-specified levels. Both these features are taken into account in \cite{Song_and_Fellouris_2016} when  all sources are observed at all times.  Thus, the present paper can be seen as a generalization of   \cite{Song_and_Fellouris_2016}  to the case that it is not possible to observe all sources at all times. However, instead of demanding  that the number of sampled sources  per  time instance be fixed, as in  \cite{Zigangirov_1966, Klimko1975OptimalSS,Dragalin_1996,oddball_2018,Cohen2015active,huang2017active,Cohen2019nonlinearcost,Tsopela_2019,Cohen2020composite},   we only require that  \textit{the ratio of the  expected number of observations over the expected time until stopping}  not exceed a  user-specified level.   This  leads to a  more general formulation for sequential anomaly detection compared to those in \cite{Zigangirov_1966, Klimko1975OptimalSS,Dragalin_1996,oddball_2018,Cohen2015active,huang2017active,Cohen2019nonlinearcost,Tsopela_2019,Cohen2020composite},  which at the same time is \textit{not}  a special case of the   sequential multi-hypothesis  testing problem with controlled sensing in \cite{chernoff1959,albert1961,Bessler1960_I,Bessler1960_II,Kiefer_Sacks_1963,Keener_1984,Lalley_Lorden_1986,nitinawarat_controlled_2013,naghshvar2013active, nitinawarat_controlled_2015,Aditya_2021}.  Thus, while existing policies  in the literature are  applicable to the proposed setup, this is not the case for the existing universal lower bounds.

Our first main result  on the proposed  problem is a criterion for a  policy (that employs  the  stopping and decision rules in \cite{Song_and_Fellouris_2016} and satisfies the  sampling constraint)  to achieve  the optimal expected time until stopping to a first-order asymptotic approximation as the two familywise  error probabilities go to 0. Indeed, we show that such a policy is asymptotically optimal  in the above sense,  if   it samples each source  with a certain minimum long-run frequency that  depends  on the source itself and  the true subset of anomalous sources. Our second main result  is that the latter condition  holds, simultaneously for every possible scenario regarding the  anomalous sources,  by  policies that  sample each source at each time instance  with a  probability that is not smaller than the  above minimum  long-run frequency that corresponds to the current estimate of the subset of anomalous  sources. This implies  the asymptotic optimality of the \textit{unmodified} policy in \cite{Bessler1960_I}, as well as of a much simpler policy according to which the sources are sampled at each time instance conditionally independent of one another given the current estimate of the subset of anomalous  sources. Indeed, the implementation of the latter policy, unlike that in \cite{Bessler1960_I}, involves minimal  computational and  storage requirements under any setup of the problem.  Moreover, we present simulation results  that suggest that this computational simplicity does not come at the price of  performance deterioration (relative to the policy in \cite{Bessler1960_I}).

Finally, to illustrate the gains of asymptotic optimality, we consider the straightforward  policy in which  the sources are sampled in tandem.  We compute  its  asymptotic relative efficiency  and we show that it  is  asymptotically optimal only in a very special setup of the problem. Moreover, our simulation results  suggest that, apart from this special setup,  its actual performance loss   relative to the above asymptotically optimal rules  is (much) larger than the one implied by its asymptotic relative efficiency  when the target error probabilities are not (very) small.

The remainder of the paper is organized as follows: in Section \ref{sec:formul} we formulate  the proposed problem.  In Section \ref{sec:policies} we present a family of policies  that satisfy the error constraints, and we introduce an auxiliary consistency property.In Section \ref{sec:probabilistic} we introduce a family of sampling rules on which we focus   in this paper.  In  Section \ref{sec:main} we  present the asymptotic optimality theory of this work. In Section  \ref{sec:other} we  discuss alternative sampling approaches in the literature, and   in Section \ref{sec:simulations} we present the results of our simulation studies. In Section \ref{sec:conclusions} we conclude and discuss  potential generalizations, as well as  directions for further research. The proofs of all  main results   are presented in  Appendices \ref{app:consistency_criterion},  \ref{app:lower_bound},   \ref{app:upper_bound}, whereas in Appendix  \ref{appen:A} we state and prove two supporting  lemmas.

We end this section with some notations we use throughout the paper. We use  $:= $ to indicate the definition of a new quantity and $\equiv$ to indicate a duplication of notation. We set 
$\bN := \{1, 2 \ldots, \}$ and  $[n] := \{1, \ldots, n\}$ for  $n \in \bN$, we denote by $A^c$ the complement, by  $|A|$ the size and by $2^A$ the powerset of a set $A$,  by $\lfloor  a  \rfloor$ the floor and by $\lceil a  \rceil$ the ceiling of a  positive number $a$, and by $\mathbf{1}$ the indicator of an event. 
We write $x\sim y$ when  $\lim (x/y)=1$,  $x\gtrsim y$ when  $\liminf (x/y) \geq 1$, and $x\lesssim y$ when  $\limsup (x/y) \leq 1$, where the limit is taken in some sense that will be specified. Moreover,   \textit{iid} stands for  independent and identically distributed, and we say that a sequence of positive numbers $(a_n)$ is \textit{summable} if $\sum_{n=1}^\infty a_n< \infty$ and \textit{exponentially decaying} if there are real numbers  $c,d>0$ such that $a_n \leq c \exp\{-d\, n \}$ for every $n \in \bN$. A property that we use   in our proofs is that  if $(a_n)$  is exponentially decaying, so is  the sequence $(\sum_{m \geq \zeta n} a_m)$, for any  $\zeta \in (0,1]$.

\section{Problem formulation}\label{sec:formul}
Let $(\bS,  \cS)$ be an arbitrary  measurable space and  let   $(\Omega, \cF, \Pro)$  be  a probability space that hosts   $M$ independent  sequences of iid $\bS$-valued random elements,  $\{X_i(n): n \in \bN\},   i \in [M]$, which are  generated by $M$ distinct data sources, as well as  an independent sequence   of iid random vectors, $\{Z(n): n =0,1, \ldots \}$,  to be used for randomization purposes. Specifically, each $Z(n) := (Z_{0}(n),Z_1(n), \ldots, Z_M(n))$ is a vector of independent random variables,  uniform in $(0,1)$, and each  $X_i(n)$ has density $f_i$, with respect to some $\sigma$-finite measure $\nu_i$,  that is equal to either $f_{1i}$ or $f_{0i}$. For every $i \in [M]$ we say that  source $i$ is  ``anomalous'' if  $f_i=f_{1i}$ and we   assume that the Kullback-Leibler divergences of   $f_{1i}$ and  $f_{0i}$    are positive and finite, i.e.,  
\begin{align}\label{KL}
	\begin{split}
		I_{i} &:= \int_{\bS} \log(f_{1i} / f_{0i}) \, f_{1i} \, d \nu_i \in (0, \infty), \\
		J_{i} &:= \int_{\bS} \log(f_{0i} / f_{1i}) \, f_{0i} \,  d \nu_i \in (0, \infty).
	\end{split}
\end{align}

We  assume that it is known \textit{a priori} that there are at least $\ell$ and at most $u$ anomalous sources, where $\ell$ and $u$ are given, user-specified integers such that  $0 \leq \ell \leq u \leq M$, with the understanding that if $\ell=u$, then $0<\ell<M$. Thus,   the  family of all possible subsets of anomalous sources is $\cP_{\ell, u}:=  \{A \subseteq [M]: \ell \leq |A| \leq u\}$. In what follows, we denote by $\Pro_A$  the underlying probability measure and by $\Exp_A$ the corresponding expectation when the subset of anomalous sources is  $A \in \cP_{\ell, u}$, and  we simply write  $\Pro$ and $\Exp$ whenever the identity of the subset of anomalous sources is not relevant. 

The  problem we consider in this work is the identification  of  the anomalous sources, if any,   on the basis of  sequentially acquired observations from all sources, when however it is not possible  to observe all of them at every sampling instance. Specifically,  we have to specify  a random  time  $T$, at which sampling is terminated, and  two random  sequences,  $R := \{R(n), n  \geq 1\}$ and  $\Delta := \{\Delta(n), n  \in \bN\}$, so that
$R(n) \subseteq [M]$ represents  the subset of sources that are sampled  at time   $n$ when  $n \leq T$, and   $\Delta(n) \equiv \Delta_n \in \cP_{\ell, u}$ represents  the subset of sources that are  identified as anomalous when $T=n$. The decisions whether to stop or not at each time instance,  which sources to sample next in the latter case,  and  which ones   to identify as anomalous in the former,  must be based  on the already available information. Thus, we say that 
$R$ is a \textit{sampling rule}  if  $R(n)$ is $\cF_{n-1}^R$--measurable  for every $n \in \bN$,  where 
\begin{align} \label{filtration}
	\begin{split}
		\cF_n^R &:=\sigma\left( \cF_{n-1}^R,\, Z(n),\, \{X_i(n):  i \in R(n)\} \right),\quad  n \in \bN, \\
		\cF_0^R &:=\sigma(Z(0)).
	\end{split}
\end{align}  
Moreover, we say  that the  triplet  $(R, T, \Delta)$ is  a \textit{policy} if $R$ is a sampling rule,  $\{T=n\} \in \cF^R_{n}$  and  $\Delta_{n}$ is $\mathcal{F}^R_{n}-$measurable  for every $n \in \bN$,  in which case we  refer to  $T$ as a \textit{stopping rule}  and to $\Delta$ as  a \textit{decision rule}. For any sampling rule $R$,  we denote by  $R_i(n)$ the indicator of whether source $i$ is sampled at time $n$, i.e., $R_i(n) := \mathbf{1}\{ i \in R(n) \}$, and by  $N_i^R(n)$ (resp. $\pi_{i}^R(n)$)  the number (resp. proportion) of times source $i$ is sampled in the first $n $ time instances,  i.e.,
\begin{align*}
	N^R_i(n) &:= \sum_{m=1}^n  R_i(m),  \quad \ \pi_i^R(n) :=  N^R_i(n)/n.
\end{align*}

We say that  a policy  $(R, T,\Delta)$ belongs to class $\cC(\alpha,\beta, \ell, u,K)$ if its probabilities  of  at least one \textit{false alarm} and at least one \textit{missed detection} upon stopping
do not exceed   $\alpha$    and  $\beta$ respectively, i.e., 
\begin{align} \label{err_const}
	\begin{split}
		\Pro_{A} \left( T <\infty, \, \Delta_{T} \setminus A \neq \emptyset \right ) \leq\alpha  \quad \forall  \, A \in \mathcal{P}_{\ell,  u}, \\
		\Pro_{A} \left( T<\infty, \, A \setminus \Delta_{T} \neq \emptyset  \right ) \leq \beta   \quad \forall   \,  A \in \mathcal{P}_{\ell,  u},
	\end{split}
\end{align}
where $\alpha, \beta$ are  user-specified  numbers in $(0,1)$, and the ratio of its expected total number of observations  over  its expected time until stopping does not exceed $K$, i.e., 
\begin{equation} \label{samp_const}
	\sum_{i=1}^M \Exp\left[  N_i^R(T) \right]  \leq K \; \Exp[T] ,
\end{equation} 
where $K$ is a user-specified, real  number in $(0,M]$. Note that, in view of the  identity 
\begin{align*}
	\sum_{i=1}^M  \Exp \left[  N_i^R(T) \right]  &= 
	\Exp \left[   \sum_{n=1}^T   \Exp \left[ |R(n)|  \; | \; \cF_{n-1}^R \right]  \right],
\end{align*}
constraint \eqref{samp_const} is  clearly satisfied when
\begin{equation} \label{samp_const_0}
	\sup_{n\leq T} \,  \Exp \left[ |R(n)|  \; | \; \cF_{n-1}^R \right] \leq  K,  
\end{equation} 
This is the case, for example, when at most  $\lfloor K \rfloor$ sources are sampled at each time instance up to stopping, i.e., when  $|R(n)| \leq \lfloor K \rfloor$ for every $n \leq T$.  

Our main goal in  this work is, for any given  $\ell, u, K$,  to obtain  policies  that attain  the smallest possible expected time until stopping, 
\begin{equation} \label{J}
	\mathcal{J}_A(\alpha, \beta,\ell, u, K) :=  \inf_{(R,T,\Delta) \in  \cC(\alpha,\beta, \ell, u,K) } \; \Exp_{A}[T],
\end{equation} 
\textit{simultaneously under every}  $A \in \mathcal{P}_{\ell, u}$,  to a first-order asymptotic approximation as   $\alpha$ and $\beta$ go to 0. Specifically,  when $\ell=u$, we allow
$\alpha$ and $\beta$ to go to 0  at arbitrary rates, but when   $\ell<u$, we assume that 
\begin{align} \label{r}
	&\exists \; r \in (0, \infty): |\log \alpha| \sim r \,  |\log \beta|.
\end{align}

\section{A family of policies} \label{sec:policies}
In this section we introduce   the statistics that we use in this work,  a family of policies that satisfy the error constraint \eqref{err_const}, as well as  an auxiliary  consistency property.

\subsection{Log-likelihood ratio statistics}
Let  $A,C \in \cP_{\ell, u}$ and
$n \in \bN$.  We denote by $\Lambda^R_{A,C} (n)$
the  log-likelihood ratio  of $\Pro_{A}$ versus $\Pro_{C}$  based on  the first $n$ time instances when the sampling rule is $R$, i.e., 
\begin{align}
	\Lambda^R_{A,C} (n) :=  \log  \frac{d\Pro_{A}}{d\Pro_{{C}}}  \left(\mathcal{F}^R_n \right),
\end{align}
and we observe that it admits the following recursion:
\begin{align}\label{LLR_global}
	\Lambda^R_{A,C} (n) &=  \Lambda^R_{A,C} (n-1)+ 
	\sum_{i \in A\setminus {C}} g_i(X_i(n)) \, R_i(n) - \sum_{j \in {C}\setminus A} g_j(X_j(n)) \, R_j(n),
\end{align}
where $ \Lambda^R_{A,C} (0) :=0$ and  $g_i :=  \log \left( f_{1i} / f_{0i}  \right), \,  i \in [M]$. Indeed, this follows by $\eqref{filtration}$ and the fact that 
$R(n)$ is $\cF_{n-1}^R$-measurable,   $X_i(n)$ is  independent of $\cF_{n-1}^R$ and its  density under $\Pro_A$ is    $f_{1i}$ if $i \in A$ and  $f_{0i}$ if $i \notin A$, and $Z(n)$ is  independent of  both $\cF_{n-1}^R$ and $\{X_i(n): i \in [M]\}$ and has the same density under $\Pro_A$ and $\Pro_C$.

{For any  $i, j \in [M]$ we write \begin{align*}
		\Lambda^R_{A,C}(n) &\equiv  \Lambda^R_{ij}(n) \quad \text{when } \quad  A=\{i\}, C=\{j\} ,\\
		\Lambda^R_{A,C}(n) &\equiv  \Lambda^R_{i}(n) \quad \text{when } \quad A=\{i\}, C=\emptyset,
	\end{align*}
	and we observe that the above recursion implies that   
	\begin{align}
		\Lambda^R_{i}(n) &= \sum_{m=1}^n   g_i(X_i(m)) \, R_i(m),  \label{LLR}\\
		\Lambda^R_{A,C} (n) &= \sum_{i \in A\setminus {C}} \Lambda^R_{i}(n) - \sum_{j \in {C}\setminus A}  \Lambda^R_{j}(n). \label{LLR_global_repre}
	\end{align}
	In particular, for  any $i,j \in [M]$ we have 
	$$
	\Lambda^R_{ij} (n) = \Lambda^R_{i}(n) -  \Lambda^R_{j}(n).$$
	
	In what follows, we refer to $\Lambda^R_{i}(n)$ as the \textit{local log-likelihood ratio} (LLR)  of source $i$ at time $n$.  We  introduce the order statistics of the LLRs at time $n$, 
	\begin{equation*}
		\Lambda^R_{(1)}(n) \geq \ldots \geq  \Lambda^R_{(M)}(n),
	\end{equation*}
	and  we denote by $w^R_{i}(n), i \in [M]$ the corresponding indices, i.e.,
	\begin{equation*}
		\Lambda^R_{(i)}(n):= \Lambda^R _{w^R_{i}(n)}(n), \quad  i \in [M].
	\end{equation*}
	Moreover, we denote by $p^R(n)$  the number of positive LLRs at time $n$, i.e., 
	$$p^R(n)\ :=  \sum_{i=1}^M \mathbf{1}  \{\Lambda^R_i(n) > 0\}, $$ 
	and we also set 
	\begin{equation*}
		\Lambda^R_{(0)}(n):= +\infty , \qquad   \Lambda^R_{(M+1)}(n) := -\infty.  
	\end{equation*}

	\subsection{Stopping and decision rules}\label{subsec:stop_deci}
	We next show that for any sampling rule $R$   there is a stopping rule, which we will denote by $T^R$,  and a decision rule, which we will denote by $\Delta^R$,   such that the policy $(R,  T^R, \Delta^R)$ satisfies the error constraint \eqref{err_const}.  Their
	forms depend on whether the number of anomalous sources is known in advance   or not, i.e., on whether $\ell=u$ or $\ell<u$.  Specifically, when $\ell=u$, we  stop as soon as the  $\ell^{th}$ largest LLR  exceeds the next one  by $c>0$, i.e., 
	\begin{align} \label{gap}
		T^{R} &:=\inf \left\{ n \in \bN: \; \Lambda^R_{(\ell)}(n)-\Lambda^R_{(\ell+1)}(n) \geq c \right\},
	\end{align}  
	and we   identify as anomalous  the sources with the $\ell$ largest LLRs, i.e., 
	\begin{align}  \label{gap_decision rule}
		\Delta^R_n &:= \left\{w^R_1(n),\ldots, w^R_\ell(n) \right\}, \quad n \in \mathbb{N}. 
	\end{align}  
	
	When the number of anomalous sources is completely unknown ($\ell=0$ and $u=M$),  we  stop as soon as the value of every LLR is outside  $(-a, b)$ for some  $a,b>0$, i.e., 
	\begin{align} \label{intersection}
		T^{R} &:= \inf \left\{n \in \bN: \;  \Lambda^R_{i}(n) \notin (-a,b)  \quad \text{for all}  \quad  i \in [M] \right\},
	\end{align}
	and we  identify  as  anomalous  the sources with positive   LLRs, i.e.,  
	\begin{align}  
		\Delta^{R}_{n} &:=\left\{ i \in [M]:  \; \Lambda^R_{i}(n)>0 \right\}, \quad n \in \mathbb{N}. 
	\end{align}
	When $\ell<u$, in general, we combine the stopping rules of  the two previous cases and we set
	\begin{align} \label{gap_intersection}
		\begin{split}
			T^{R} := \inf  \{n \in \bN:  \quad &   \text{either} \quad  \Lambda^R_{(\ell+1)}(n){\leq}-a  
			\quad  \& \quad  \Lambda^R_{(\ell)}(n)-\Lambda^R_{(\ell+1)}(n) \geq c,  \\
			\quad & \; \text{or} \quad   \quad  \ell \leq p^{R}(n) \leq u \qquad  \& \qquad   \Lambda^R_{i}(n) \notin (-a,b)  \quad \forall \; 
			i \in [M],   \\
			\quad  & \; \text{or} \quad  \quad \;  \Lambda^R _{(u)}(n) \geq b \qquad  \& \quad  \Lambda^R_{(u)}(n)-\Lambda^R_{(u+1)}(n) \geq d \},  
		\end{split}
	\end{align}
	where $a,b,c,d>0$,  and we use the following decision rule:
	\begin{align} \label{gi_decision rule}
		\Delta^{R}_{n} &:=\left\{ w^R_{i}(n): \; i=1, \ldots, (p^R(n) \vee \ell) \wedge u \right\}, \quad n \in \mathbb{N}.
	\end{align}
	That is,  we identify as anomalous the sources with positive LLRs as long as their number is between $\ell$ and $u$. If this number is larger than $u$ (resp. smaller than $\ell$), then we declare as anomalous the sources  with the $u$ (resp. $\ell$) largest LLRs. \\

	\begin{proposition} 
		Let  $R$ be an arbitrary sampling rule. 
		\begin{itemize}
			\item When  $\ell=u$,  $(R, T^{R}, \Delta^{R})$  satisfies the error constraint \eqref{err_const}   if  
			\begin{align} \label{thresholds_gap}
				c=  |\log(\alpha \wedge \beta)|+\log(\ell(M-\ell)).
			\end{align} 
			\item When $\ell<u$,    $(R, T^{R}, \Delta^{R})$  satisfies the error constraint \eqref{err_const}   if  
			\begin{align} \label{thresholds_gi}
				\begin{split}
					a&=|\log\beta|+\log M, \quad c=|\log \alpha|+\log((M-\ell)M), \\
					b&=|\log \alpha|+\log M, \quad d=|\log \beta|+\log(uM).
				\end{split}
			\end{align}
		\end{itemize}
	\end{proposition}
	
	\begin{IEEEproof}
		When all sources are sampled at all times, this is shown in \cite[Theorems 3.1, 3.2]{Song_and_Fellouris_2016}. The same proof applies when $K<M$ for any sampling rule  $R$. \\
	\end{IEEEproof}

	In view of the above result,  in what follows we assume that the  thresholds  in $T^R$ are selected according to  \eqref{thresholds_gap} when $\ell=u$  and according to  \eqref{thresholds_gi} when $\ell<u$.  While this is a rather conservative choice, it will be sufficient for obtaining asymptotically optimal policies. For this reason, in what follows we  say that a sampling rule, $R$, that satisfies the sampling constraint \eqref{samp_const} with $T=T^R$, is \textit{asymptotically optimal under $\Pro_A$}, for some $A \in \cP_{\ell, u}$,   if 
	$$\Exp_{A}\left[T^R \right] \sim \cJ_A(\alpha, \beta,\ell, u, K)$$ 
	as $\alpha$ and $\beta$ go to $0$ at arbitrary rates when $\ell=u$ and so that  \eqref{r} holds  when $\ell<u$.  
	We simply say that   $R$ is \textit{asymptotically optimal} if it is
	asymptotically optimal under $\Pro_A$ for every $A \in \cP_{\ell,u}$.  We next introduce a weaker property, which will be useful for establishing asymptotic optimality.

	\subsection{Exponential consistency}   For any sampling rule    $R$
	and any subset $A \in \cP_{\ell, u}$ we denote by $\sigma_A^R$ the random time starting from which the sources in $A$ are the ones  estimated as  anomalous by $\Delta^R$, i.e.,
	\begin{equation}\label{sigmaA}
		\sigma^R_{A} :=\inf\left\{n \in \bN: \Delta^{R}_{m} =A  \quad \text{for all} \; \;   m \geq n\right\},
	\end{equation} 
	and we say that  $R$ is \textit{exponentially consistent under $\Pro_A$} if $\Pro_A(\sigma_A^R >  n)$ is an exponentially decaying  sequence. We simply say that   $R$ is \textit{exponentially consistent} if it is
	exponentially consistent under $\Pro_A$ for every $A \in \cP_{\ell,u}$. The following theorem  states sufficient conditions for exponential consistency under $\Pro_A$. \\

	\begin{theorem}\label{propo_1}
		Let $A \in \cP_{\ell,u}$ and let $R$ be an arbitrary sampling rule. 
		\begin{itemize} 
			\item   When  $\ell<u$,  $R$ is  exponentially consistent under $\Pro_A$ if  there exists a  $ \rho>0$ such that  
			$\Pro_{A} \left( \pi^R_i(n) <  \rho  \right)$ is exponentially decaying  for every $i \in A$ if   $ |A|>\ell$
			and for  every $i \notin A$ if $ |A|<u$.
			\item   When  $\ell=u$,  $R$ is  exponentially consistent under $\Pro_A$ if  there exists a  $ \rho>0$ such that  
			$\Pro_{A} \left( \pi^R_i(n) <  \rho  \right)$ is exponentially decaying   either for every $i \in A$ or for every $i \notin A$.
		\end{itemize}
		
	\end{theorem}
	
	\begin{IEEEproof} 
		Appendix \ref{app:consistency_criterion}.\\
	\end{IEEEproof}

	\noindent {\underline{\textbf{Remark:}} Theorem \ref{propo_1} reveals that when  $|A|=\ell>0$ or $|A|=u<M$, it is possible to have exponentially consistency under $\Pro_A$  without sampling at all certain sources. Specifically, when $|A|=\ell<u$ (resp. $|A|=u>\ell$) it is not necessary to sample any source in $A$ (resp. $A^c$). On the other hand, when $|A|=\ell=u$, it suffices to  sample either all sources in $A$ or all of them in $A^c$.\\
	}

	\section{Probabilistic sampling rules}  \label{sec:probabilistic}
	In this section we introduce a family of sampling rules and we show how to design them in order to satisfy the sampling constraint \eqref{samp_const_0} and  in order to be exponentially consistent. Thus, we say that a sampling rule $R$ is  \textit{probabilistic} if there exists a  function $q^R : 2^{[M]} \times \cP_{\ell, u} \to [0,1]$
	such that    for every  $n  \in \bN$, $D \in \cP_{\ell,u}$,  and  $B  \subseteq [M]$ we have   
	\begin{align} \label{q}
		q^R\left(B;D \right) := \Pro \left( R(n+1)=B \,  |  \, \cF^R_{n} , \Delta^{R}_{n}=D\right),
	\end{align} 
	i.e., $q^R\left(B;D \right)$ is  the probability that $B$  is the subset of sampled sources   when $D$ is the  currently estimated  subset of anomalous sources.
	If $R$ is a probabilistic sampling rule,  then for every  $i \in [M]$ and $D \in \cP_{\ell,u}$ we set
	\begin{align} \label{c}
		\begin{split}
			c^R_{i} \left(D \right) &:=\Pro \left( R_i(n+1)=1 \, |  \, \cF^R_{n}, \Delta^R_{n}=D \right)\\
			&=  \sum_{B \subseteq [M]: \, i \in B}  q^R\left(B;D \right),
		\end{split}
	\end{align}
	i.e.,  $c^R_i(D)$ is the probability that  source $i$ is sampled when $D$ is the currently  estimated  subset of anomalous sources.

	We refer to a probabilistic sampling rule $R$  as \textit{Bernoulli}  if  for every $D \in \cP_{\ell, u}$ and 
	$B  \subseteq [M]$ we have 
	\begin{align}  \label{bernoulli}
		q^R(B; D)= \prod_{i \in B} c^R_i(D)  \prod_{j \notin B} \left(1-c^R_j(D)\right), 
	\end{align}
	i.e., if the  sources are sampled at each time instance  conditionally independent of one another given the currently estimated subset of anomalous sources. Indeed, such a sampling rule  admits a representation of the form 
	\begin{align}\label{bern_1} 
		R_i(n+1) &= \mathbf{1} \left\{Z_i(n)\leq c^R_i\left(\Delta_{n}^R \right) \right\}, \quad i \in [M], n \in \bN, 
	\end{align} 
	where  $Z_1(n), \ldots, Z_M(n)$ are iid and uniform in $(0,1)$,
	thus, its   implementation at each time instance requires  the realization of $M$ Bernoulli random variables.\\

	The following proposition provides  a sufficient condition for a probabilistic sampling rule to satisfy the sampling constraint  \eqref{samp_const_0}, and consequently \eqref{samp_const},  for any $\{\cF_n^R\}$-stopping time $T$. 
	\\
	
	\begin{proposition} \label{prop:prob_samp_constraint}
		If $R$ is a  probabilistic sampling rule such that 
		\begin{align} \label{sum_less_K} 
			\sum_{i=1}^M  & c^R_i(D) \leq K \quad \text{for every } \; D \in \cP_{\ell, u},
		\end{align}
		then \eqref{samp_const_0} holds for any $\{\cF_n^R\}$-stopping time $T$. 
	\end{proposition}
	
	\begin{IEEEproof}
		For any   $n \in \bN$ and any probabilistic sampling rule $R$ we have 
		$$		 \Exp \left[ |R(n)|  \; | \; \cF_{n-1}^R \right] = \sum_{i=1}^M \Pro \left(R_i(n)=1  \; | \; \cF_{n-1}^R \right)= \sum_{i=1}^M  c^R_i(\Delta_{n}^R).$$
		As a result, if  \eqref{sum_less_K}  is satisfied, then   \eqref{samp_const_0} holds for any $\{\cF_n^R\}$-stopping time, $T$.\\
	\end{IEEEproof}

	Finally, we establish sufficient conditions for the exponentially consistency of a probabilistic sampling rule.

	\begin{theorem}\label{theo:prob_cons} Let  $R$ be a  probabilistic sampling rule. 
		\begin{itemize} 
			\item   When $\ell<u$,  $R$ is  exponentially  consistent  if,  for every $D  \in \cP_{\ell,u}$,   $c_i^R(D)$ is positive  
			for every $i \in D$ if $|D|>\ell$ and  every $i \notin D$  if    $|D|<u$.
			
			\item   When $\ell=u$,  $R$ is  exponentially  consistent  if,  for every $D  \in \cP_{\ell,u}$,   $c_i^R(D)$ is positive   either 
			for every $i \in D$ or for every $i\notin D$.
		\end{itemize}
	\end{theorem}
	
	\begin{IEEEproof}  The proof consists in showing that the sufficient conditions of  Theorem \ref{propo_1} are satisfied for every $A \in \cP_{\ell,u}$, and is presented in Appendix \ref{app:consistency_criterion}.  \\
	\end{IEEEproof}

	\noindent { \underline{\textbf{Remark:}}  Theorem \ref{theo:prob_cons}  implies that when $\ell<u$, a probabilistic sampling rule is exponentially consistent  if, whenever the number of  estimated anomalous  sources is larger than $\ell$ (resp. smaller than $u$), it samples with positive probability   any source that is currently estimated as anomalous  (resp. non-anomalous). When $\ell=u$, on the other hand, 
		it suffices to sample with positive probability  at any time instance  either every source that is currently estimated as anomalous  or every source that is  currently estimated as  non-anomalous.} \\
	
	\noindent { \underline{\textbf{Remark:}} Unlike the proof of Proposition \ref{prop:prob_samp_constraint}, the proof of Theorem \ref{theo:prob_cons}  is quite challenging and it is one of the main contributions of this work from a technical point of view. It relies on two Lemmas,   \ref{lem:prob_consi_lem0}  and  \ref{lem:prob_consi_lem1},
		which we state  in  Appendix \ref{app:consistency_criterion}. The proof becomes much easier if we strengthen the assumption of the theorem and assume that every source be sampled with positive probability at every time instance, i.e., if we assumed that  $c_i^R(D)>0$ for every $D \in \cP_{\ell,u}$. Indeed, in this case  Lemma \ref{lem:prob_consi_lem1} becomes redundant,  the proof simplifies considerably, and it essentially relies on ideas developed in \cite{chernoff1959}.  However,  such an assumption would limit considerably the scope of the  asymptotic optimality theory we develop in the next section. }

\section{Asymptotic optimality} \label{sec:main}

In this section we present the asymptotic optimality theory  of this work and discuss some of its implications. For this, we first  need to introduce some additional notation.
\subsection{Notation}
For $A \subseteq [M]$ with $A \neq \emptyset$ we  set
\begin{align}\label{KL_A}
	I^{*}_{A} &:= \min\limits_{i \in A}I_{i} , \quad 
	I_{A} := \frac{|A|}{\sum_{i \in A} (1/I_{i})}, \quad \hat{K}_{A} :=|A| \; \frac{I^{*}_{A}}{ I_{A}},
\end{align}
and for $A \subseteq [M]$ with  $A \neq [M]$  we set
\begin{align} \label{KL_notA}
	J^{*}_{A}:= \min\limits_{i \notin A} J_{i} \quad J_{A}:= \frac{|A^c|}{\sum_{i \notin A}(1/J_{i})}, \quad \check{K}_{A} := |A^c| \;  \frac{ J^{*}_{A} }{ J_{A}}.
\end{align} 
That is,  $I^{*}_{A}$  is the minimum and   $I_{A}$  the harmonic mean  of the Kullback-Leibler information numbers in $\{I_{i}: i \in A\}$,  whereas $J^{*}_{A}$  is the minimum and   $J_{A}$  the harmonic mean  of the Kullback-Leibler information numbers in $\{J_{i} : i \notin A\}$. Moreover,  $\hat{K}_{A}$ (resp. $\check{K}_{A}$) is
a positive real  number smaller or equal to  the size of $A$ (resp. $A^c$), with the equality attained when $I_i=I$  (resp. $J_i=J$) for every $i$ in  $A$  (resp. $A^c$). 

Moreover,  for $A \subseteq [M]$ with $0<|A|<M$
we set 
\begin{align}\label{theta_def}
	\theta_{A} &:=\;  I^{*}_{A} /  J^{*}_{A},
\end{align}
thus, $\theta_A$ is a positive real  number that is equal to 1 when  $I_i=J_i$ for every $i \in [M]$. 

In   Theorem \ref{theo:lowerbound} we also   introduce for each $A \in \cP_{\ell, u}$ two quantities, 
\begin{equation}
	x_A(r,\ell, u, K) \quad \text{and} \quad y_A(r,\ell,u, K),
\end{equation}
which will play an important role in the formulation of all results in this section.  Although their values  depend on  $r, \ell, u, K$, in  order to lighten the notation   we simply write $x_A$ and $y_A$  unless we want to emphasize their dependence on one or more of these parameters. 


\subsection{A universal asymptotic lower bound} \label{subsec:lower_bound}

\begin{theorem} \label{theo:lowerbound}
	Let  $A \in \cP_{\ell, u}$. 
	
	(i) Suppose that  $ \ell=u $. Then, 
	as  $\alpha, \beta \to 0$ we have  
	\begin{align}  \label{LB_gap}
		\mathcal{J}_A(\alpha, \beta,\ell, u, K) &\gtrsim \frac{|\log ( \alpha \wedge \beta)|}{x_A\,  I^{*}_{A}  +y_A  \,  J^{*}_{A}} ,
	\end{align}	
	where  if  $\hat{K}_A  \leq \theta_A \check{K}_A$,  
	\begin{align} \label{xy_gap1}
		\begin{split}
			x_A&:=(K/ \hat{K}_{A} ) \wedge    1 \\
			y_A &:=  \bigl( ( K-\hat{K}_{A} )^{+} /  \check{K}_{A} \bigr) \wedge 1, 
		\end{split}
	\end{align}	
	and  if  $\hat{K}_A  > \theta_A \check{K}_A$,  
	\begin{align} \label{xy_gap}
		\begin{split}
			x_A&:= 
			\bigl(  (K-\check{K}_{A} )^{+} /  \hat{K}_{A}  \bigr)  \wedge 1 
			\\
			y_A&:= (K/ \check{K}_{A} ) \wedge    1 .
		\end{split}
	\end{align}

	(ii) Suppose  that  $ \ell<u$ 
	and let   $\alpha, \beta \to 0$ so that \eqref{r} holds.

	\begin{itemize}		
		
		\item  If  $\ell <|A| < u$, then
		\begin{align}  \label{LB_int}
			\begin{split}
				\mathcal{J}_A(\alpha, \beta,\ell, u, K)&\gtrsim \frac{ |\log \alpha| }{x_A\,  I_A^*} 
				\sim \frac{|\log \beta| }{y_A\,  J_A^*}, \\ \text{where} \qquad x_A&:= \frac{K}{ \hat{K}_A+(\theta_{A}/r) \check{K}_A} \wedge (r/\theta_A) \wedge 1   \\
				y_A&:= \frac{K}{ \check{K}_A+(r/\theta_{A}) \hat{K}_A} \wedge (\theta_A/r) \wedge 1 . 
			\end{split}
		\end{align}

		\item If $|A|=\ell$,  we distinguish two cases:
		\begin{itemize}
			\item[] If   either  $\ell=0$ or   $r \leq 1$, then 
			\begin{align}\label{LB_00}
				\begin{split}
					\mathcal{J}_A(\alpha, \beta,\ell, u, K) &\gtrsim  \frac{|\log \beta|}{x_A \, I_A^*+ y_A \, J_A^*} , \quad \\
					\text{where} \qquad x_{A}&:=0 \\
					y_A &:= (K/ \check{K}_A) \wedge 1. 
				\end{split}
			\end{align}			
			If  $\ell>0$ and  $r>1$, we set $z_A:=\theta_{A}/(r-1)$ and we distinguish two further  cases:
			
			\begin{itemize}
				
				\item[] \quad  If  either $z_A  \geq 1$  \textit{or}  $K \leq \hat{K}_{A}+z_A \, \check{K}_A$, then 
				\begin{align} \label{LB_gap_int_1}
					\begin{split}
						\mathcal{J}_A(\alpha, \beta,\ell, u, K) &\gtrsim   
						\frac{|\log \beta | }{y_A\,  J_A^*} 
						\sim   \frac{|\log \alpha | }{x_A\,  I_A^*+y_A \, J_A^*},   
						\\
						\text{where} \qquad   
						x_{A}&:= \frac{K}{\hat{K}_{A}+z_A \, \check{K}_A}   \wedge  (1/z_A)\wedge  1  \\
						y_{A}&:= \frac{K}{\check{K}_{A}+(1/z_A) \, \hat{K}_A}   \wedge  z_A \wedge  1. 
					\end{split}
				\end{align}	
				
				\item[] \quad  If \,  $z_A <1$  \textit{and}   $K >\hat{K}_{A}+z_A \, \check{K}_A$, then 
				\begin{align} \label{LB_gap_int_2}
					\begin{split}
						\mathcal{J}_A(\alpha, \beta,\ell, u, K)  &\gtrsim \frac{| \log \alpha|}{x_A\,  I^{*}_{A}  +y_A  \,  J^{*}_{A}}  \\
						\; \; \text{where}  \qquad              x_{A}&:=1 \\
						y_{A}&:= \bigl( (K-\hat{K}_{A} )/  \check{K}_{A} \bigr) \wedge 1.
					\end{split}
				\end{align}

			\end{itemize}	
			
		\end{itemize}

		\item If $|A|=u$, then we distinguish again two cases:
		
		\begin{itemize}
			\item[] If  either $u=M$ or $r \geq 1$, then 
			\begin{align} \label{LB_11}
				\begin{split}
					\mathcal{J}_A(\alpha, \beta,\ell, u, K) &\gtrsim \frac{|\log \alpha|} {x_A \, I_A^*+ y_A \, J_A^*}  \\
					\text{where} \qquad  y_{A}&:=0 \\
					x_A &:= (K/ \hat{K}_A) \wedge 1.	
				\end{split} 
			\end{align}	
			If 	$u<M$ and $r <1$,  we set $w_A:= (1/\theta_A)/( 1/r -1)$ and  we distinguish two  further  cases:
			
			\begin{itemize}	
				\item[] \quad If  either $w_A\geq 1$ or $K  \leq   \check{K}_A+ w_A \hat{K}_A$, then
				\begin{align} \label{LB_gap_int_11}
					\begin{split}
						\mathcal{J}_A(\alpha, \beta,\ell, u, K)&\gtrsim \frac{|\log \alpha | }{x_A \, I_A^*}  \sim   \frac{|\log \beta|}{x_A \, I_A^*+ y_A \, J_A^*} ,  \\  \text{where} \qquad  y_A &:= 
						\frac{K}{ \check{K}_A+ w_A \hat{K}_A} \wedge (1/w_A) \wedge 1  \\
						x_A &:= 
						\frac{K}{ \hat{K}_A+ (1/w_A) \check{K}_A} \wedge w_A \wedge 1. 
					\end{split} 
				\end{align}	
				
				\item[] \quad If  $w_A<1$ \textit{and} $K >  \check{K}_A+ w_A \hat{K}_A$, then 
				\begin{align} \label{LB_gap_int_12}
					\begin{split}
						\mathcal{J}_A(\alpha, \beta,\ell, u, K) &\gtrsim \frac{|\log \beta|}{x_A \,  I_A^*+ y_A \, J_A^*} ,  \\
						\text{where} \qquad         y_{A}&:=1   \\
						x_{A} &:=\bigl( (K-\check{K}_{A}) / \hat{K}_{A} \bigr)  \wedge 1 . 
					\end{split}
				\end{align}						
				
			\end{itemize}
		\end{itemize}
		
	\end{itemize}
\end{theorem}

\begin{IEEEproof}  The proof  is presented in Appendix \ref{app:lower_bound}. It follows similar steps as the one in the full sampling case in 
	\cite[Theorem 5.1]{Song_and_Fellouris_2016}, with the difference that it uses a version of Doob's optional sampling theorem in the place of Wald's identity and requires the solution of a max-min optimization problem, in each of the cases we distinguish, to determine the denominator in the lower bound.\\
\end{IEEEproof}

\noindent {\textbf{\underline{Remark}:} By the definition of $x_A$ and $y_A$ in Theorem  \ref{theo:lowerbound} we can see that  
	\begin{itemize}
		\item they  both take values in  $[0,1]$ and  at least one of them is positive,
		\item they are both increasing as functions of $K$, 
		\item at least one of them is equal to 1 when $K=M$, 
		\item if $x_A=0$ (resp. $y_A=0$), then  $|A|=\ell$ (resp. $|A|=u$).
	\end{itemize} 
	In the next section we obtain an interpretation of the values of $x_A$ and $y_A$. 
}

\subsection{A  criterion for asymptotic optimality}
Based on the universal asymptotic lower bound of Theorem \ref{theo:lowerbound}, we next establish the asymptotic optimality under $\Pro_A$  of a sampling rule $R$ that satisfies the sampling constraint and  samples each source $i \in [M]$, when the true subset of anomalous sources is $A$, with a long-run frequency that is not smaller than
\begin{equation} \label{c_star}
	c_i^*(A) := 
	\begin{cases}
		x_{A} \; I_A^*/ I_i  \quad  \text{if}  \; \;  i \in A,\\
		y_{A} \; J_A^*/J_i  \quad  \text{if}  \; \;  i \notin A.
	\end{cases}
\end{equation}

\begin{theorem} \label{theo:AO} 
	Let  $A \in \cP_{\ell, u}$ and let    $R$ be a sampling rule that satisfies  \eqref{samp_const}  with $T=T^R$.
	If for every  $i \in [M]$ such that $c_i^*(A)>0$ the sequence  $\Pro_{A} \left(\pi^R_i(n) < \rho \right)$ is summable  for every   $\rho \in (0, c_i^*(A) )$,  then $R$ is  asymptotically optimal under $\Pro_A$.
\end{theorem}

\begin{IEEEproof} The proof consists in   establishing  asymptotic upper bounds for $\Exp_A[T^R]$,  when the thresholds are selected according to \eqref{thresholds_gap}-\eqref{thresholds_gi},  that match the universal asymptotic lower bounds of   Theorem \ref{theo:lowerbound}. The proof  is presented
	Appendix \ref{app:upper_bound}.\\
\end{IEEEproof}

\noindent
{ \textbf{\underline{Remark}:}  
	Let $A \in \cP_{\ell, u}$. By the   definition of $\{c_i^*(A): i \in [M]\}$    in  \eqref{c_star} it follows that:
	\begin{itemize}
		\item $c_i^*(A) \in [0,1]$ for every  $i \in [M]$, since   $x_A, y_A \in [0,1]$,
		\item  $x_A=0$ (resp. $y_A=0$) 
		$\Leftrightarrow$  $c_i^*(A)=0$
		for every $i$ in $A$ (resp. $A^c$),
		\item  $x_A>0$ (resp. $y_A>0$) 
		$\Leftrightarrow$  $c_i^*(A)>0$
		for every $i$ in $A$ (resp. $A^c$). 
	\end{itemize}
	Therefore, the above theorem implies that when $x_A$ (resp. $y_A$) is equal to $0$,  it is not necessary to sample any source in $A$ (resp. $A^c$)  in order to achieve asymptotic optimality  under $\Pro_A$. Moreover, recalling the definitions of $I_A$ and $J_A$ in \eqref{KL_A} and \eqref{KL_notA}, we can  see that:
	\begin{align} \label{identities}
		\begin{split}
			x_A&= \frac{I_A/ I^*_A}{|A|} \sum_{i \in A} \; c_i^*(A)  \quad \text{when} \quad A \neq \emptyset, \\
			y_A &=  \frac{J_A/ J^*_A}{|A^c|} \sum_{i \notin A} c_i^*(A)  \quad \text{when} \quad A \neq [M].
		\end{split}
	\end{align}
	These  identities reveal that $x_A$ (resp. $y_A$) is equal to the average of the minimum limiting sampling frequencies of the sources in $A$ (resp. $A^c$)  required for asymptotic optimality under $\Pro_A$, multiplied by $I_A/I_A^*$ (resp.  $J_A/J_A^*$).} \\

\noindent
{ \textbf{\underline{Remark}:}  
	By  the  definitions of
	$\hat{K}_A$ and $\check{K}_A$ in \eqref{KL_A}  and   \eqref{KL_notA}, we can see that \eqref{identities} implies 
	\begin{align} \label{identity_new}
		\sum_{i=1}^M c_i^*(A) &=x_A \hat{K}_A+ y_A\check{K}_A .
	\end{align} 
	Moreover, by a direct inspection of the values of $x_A$ and $y_A$ we have \begin{equation} \label{new_ineq}
		x_A \hat{K}_A+ y_A\check{K}_A \leq K,
	\end{equation}
	and consequently:
	\begin{align} \label{inequality_new}
		\sum_{i=1}^M c^*_i(A) \leq K.
\end{align} }

\subsection{Asymptotically optimal probabilistic sampling rules}

We next specialize  Theorem \ref{theo:AO}  to the case of probabilistic sampling rules. Specifically,  we obtain a sufficient condition for the asymptotic optimality, simultaneously under every possible scenario, of  an arbitrary  probabilistic sampling rule, $R$, in terms of the quantities $\{c_{i}^{R}(A): 
i \in [M], A \in \cP_{\ell,u} \}$,  defined in \eqref{c}. \\

\begin{theorem}\label{theo:prob_ao}  
	If  $R$ is a  probabilistic sampling rule and for every  $A \in \cP_{\ell, u}$  it   satisfies 
	\begin{align} \label{system}
		c^R_i(A) \geq   c_i^*(A),\quad \forall \;  i \in [M],
	\end{align}
	then	it   is exponentially consistent. If also $R$  satisfies \eqref{sum_less_K}, then it is asymptotically optimal under $\Pro_A$ for every  $A \in \cP_{\ell, u}$.
\end{theorem}

\begin{IEEEproof}  Exponential consistency is established by showing that  the conditions of Theorem \ref{theo:prob_cons} are satisfied, and  asymptotic optimality  by showing that  the conditions of Theorem \ref{theo:AO} are satisfied. The proof is presented in Appendix \ref{app:upper_bound}.\\
\end{IEEEproof}

In the next corollary we show that  both conditions of Theorem \ref{theo:prob_ao}  are satisfied when the equality holds in \eqref{system} for every $A \in \cP_{\ell,u}$. 


\begin{corollary}\label{coro:prob_ao}  
	If  $R$ is a probabilistic sampling rule  such  that   for every   $A \in \cP_{\ell, u}$ we have 
	\begin{align} \label{system_epsilon}
		c^R_i(A) = c^*_i(A), \quad  \forall \;  i \in [M],
	\end{align}
	then  $R$  is exponentially consistent and  asymptotically optimal under $\Pro_A$ for every  $A \in \cP_{\ell, u}$.
\end{corollary}

\begin{IEEEproof} By Theorem \ref{theo:prob_ao}  it clearly suffices to show that \eqref{sum_less_K} holds, which follows directly by  \eqref{inequality_new}. \\
\end{IEEEproof}

While  \eqref{system_epsilon} suffices for the asymptotic optimality of a probabilistic sampling rule under $\Pro_A$, it is not always necessary. In the following proposition we characterize the setups for which the necessity holds.  \\

\begin{proposition} \label{propo:prob_ao} 
	Let  $A \in \cP_{\ell,u}$. Then:
	\begin{center} \eqref{system_epsilon} holds for any  probabilistic sampling rule $R$  that satisfies   \eqref{sum_less_K} and \eqref{system}
	\end{center}
	\begin{center}  $\Leftrightarrow$  $x_A \hat{K}_A + y_A \check{K}_A  =K$.
	\end{center}
\end{proposition}

\begin{IEEEproof}  For any probabilistic sampling rule  $R$  that satisfies \eqref{sum_less_K}  and  \eqref{system}  we have 
	\begin{equation}\label{obser1} 
		K\geq \sum_{i=1}^M c^R_i(A) \geq \sum_{i=1}^{M}c_{i}^{*}(A) =x_A \hat{K}_A + y_A \check{K}_A,
	\end{equation}
	where the equality follows by  \eqref{identity_new}.

	$(\Leftarrow)$  If  $K=  x_A \hat{K}_A + y_A \check{K}_A$, then  by  \eqref{obser1} we obtain 
	\begin{align*}
		\sum_{i=1}^M c^R_i(A)= \sum_{i=1}^M c^*_i(A).
	\end{align*}
	In view of \eqref{system},  this proves that $ c^R_i(A) =   c_i^*(A)$ for every $i \in [M]$.

	$(\Rightarrow)$  We  argue by contradiction and assume that   $x_A \hat{K}_A + y_A \check{K}_A = K$ does not hold.  By \eqref{identity_new} and   \eqref{new_ineq} it then follows that   $\sum_{i=1}^M c^*_i(A) <K$.  This implies that  there  is a probabilistic sampling rule $R$  that satisfies  \eqref{system}  \textit{with strict  inequality for at least one $i \in [M]$}, 
	and also  satisfies \eqref{sum_less_K}. Thus, we have reached  a contradiction, which completes the proof. 
\end{IEEEproof}


\subsection{The optimal performance under full sampling}
Corollary \ref{coro:prob_ao} implies that the asymptotic lower bound  in  Theorem \ref{theo:lowerbound} is always achieved, and as  a result
it is a  first-order asymptotic approximation  to the optimal expected time until stopping. By  this approximation it follows that   for any $A \in \cP_{\ell, u}$ we have 
\begin{align}  \label{Q_equiv}
	\cJ_A(\alpha, \beta, \ell, u, K) \sim \cJ_A(\alpha, \beta, \ell, u, M)
	\; &\Leftrightarrow  \;    K   \geq Q_A,
\end{align}
where $Q_A \equiv Q_A(r, \ell, u)$ is defined as follows:
\begin{align}  \label{Q1}
	Q_A &:=  x_{A}(r, \ell,u, M) \hat{K}_{A} +y_{A}(r,\ell,u,M) \check{K}_{A}.
\end{align} 
Moreover, by an inspection of the values of $x_A$ and $y_A$ in Theorem \ref{theo:AO} it follows that  \begin{align}  \label{Q2}
	\begin{split}
		Q_A &= \min\{K \in (0,M]:  \quad  x_{A}(r, \ell,u, K)=  x_{A}(r, \ell,u, M) \\
		& \qquad \qquad \qquad \qquad \quad \&  \quad    y_{A}(r, \ell,u, K)=  y_{A}(r, \ell,u, M) \}.
	\end{split}
\end{align}
The equivalence in  \eqref{Q_equiv}  implies that  if  $Q_A <M$ and $K \in [Q_A, M)$, then the optimal  expected time until stopping under $\Pro_A$   \textit{in the case of full sampling}  can be  achieved, to a first-order asymptotic approximation as $\alpha, \beta \to 0$,   \textit{without} sampling all sources at all times. \\

\begin{corollary}
	If   $Q_A <M$ and $K \in [Q_A, M)$  for some $A \in \cP_{\ell, u}$, then there is a probabilistic sampling rule $R$ such that  
	$(R, T^R, \Delta^R) \in \cC(\alpha, \beta, \ell, u, K)$ and 
	$\Exp_A[T^R] \sim \cJ_A(\alpha, \beta, \ell, u, M),$ while  $c_i^R(A)<1$ for some $i \in [M]$. 
\end{corollary}

\begin{IEEEproof}
	This follows by Corollary \ref{coro:prob_ao} and \eqref{Q_equiv}. \\
\end{IEEEproof}

The next proposition,  in conjunction with Proposition \ref{propo:prob_ao},  shows that $Q_A$ can also be used to characterize  the  setups for which the equality must hold in \eqref{system}. In particular, it shows that   this is always the case when $K\leq 1$.

\begin{proposition}\label{propo:prob_ao2} 
	For every $A \in \cP_{\ell, u}$ we have  $Q_A \geq 1$ and 
	$$x_A \hat{K}_A + y_A \check{K}_A  =K \; \Leftrightarrow  \; K\leq Q_A.$$
\end{proposition}

\begin{IEEEproof}  By  the definition of 
	$\hat{K}_A, \check{K}_A$ in  \eqref{KL_A}, \eqref{KL_notA} it follows that  
	\begin{align*}
		\hat{K}_A &= \sum_{i \in A} (I_A^*/I_i) \geq 1 \quad \text{for} \; A \neq \emptyset \\
		\check{K}_A &= \sum_{i \notin A} (J_A^*/J_i) \geq 1 \quad \text{for} \; A \neq [M].
	\end{align*}
	Since also at least one of  $x_A(r, \ell, u, M)$  and $y_A(r, \ell, u, M)$ is equal to $1$, by the definition of $Q_A$ in  \eqref{Q1} we conclude  that $ Q_A \geq 1$.  It remains to prove the equivalence. $(\Rightarrow)$ It suffices to observe that since   $x_A$, $y_A$     are both  increasing in $K$,  by the definition of  $Q_A$ in  \eqref{Q1} we have 
	$x_A\hat{K}_A + y_A\check{K}_A  \leq Q_A$.   $(\Leftarrow)$ Follows by  direct verification.
\end{IEEEproof}

\subsection{The Chernoff sampling rule} \label{subsec:chernoff}
When $K$ is an integer,  we will refer to  a probabilistic sampling rule  as  \textit{Chernoff} if  it  satisfies the conditions of Theorem \ref{theo:prob_ao} and  samples \textit{exactly $K$ sources per  time instance}. Indeed, such a  sampling rule is  implied from  \cite{chernoff1959},\cite{Bessler1960_I}, \cite{Bessler1960_II} when these works are applied to the sequential anomaly detection problem (with a fixed number of sampled sources per time instance). In fact,   if  the class $\cC(\alpha, \beta, \ell, u, K)$ is restricted to policies  that sample \textit{exactly} $K$ sources per time instance, the asymptotic optimality of this rule under $\Pro_A$   is implied by the general results in  \cite{Bessler1960_I}, \textit{as long as  $x_A$ and $y_A$ are both positive}. However, this is not always the case.  Indeed, in the simplest formulation of the sequential anomaly detection problem, where it is known a priori that there is only one anomaly ($l=u=1$),  only one source can be sampled per time instance $(K=1)$, and the sources are homogeneous, i.e.,  $f_{0i}=f_0$, $f_{1i}=f_1$ for every $i \in [M]$, then one of $x_A$ and  $y_A$ is 0 for any  $A \in \cP_{\ell,u}$.  In this setup, the asymptotic optimality of a Chernoff rule  was shown in 
\cite[Appendix A]{Bessler1960_II} if also  $f_0$, $f_1$ are both Bernoulli pmf's.  Our  results in  this section remove all these restrictions and  establish the  asymptotic optimality of the Chernoff rule for any values of  $\ell, u, r, K$, and without artificially modifying it  at a subsequence of sampling instances, as in \cite{nitinawarat_controlled_2013}.   

On the other hand, to implement a Chernoff rule one needs to determine a function $q^R$ that   satisfies simultaneously the conditions of Theorem \ref{theo:prob_ao} and  also 
\begin{align}  \label{chernoff}
	q^R\left(B; D \right) &=0 \quad   \text{for all}  \; D \in \cP_{\ell, u} \quad  \text{and}  \quad  B \subseteq [M] \; \text{with} \;  |B| \neq  K,
\end{align}
which can be a  computationally demanding task unless the problem has a special  structure or $K=1$. This should be  compared with the corresponding asymptotically optimal \textit{Bernoulli} sampling rule, whose implementation does not  require essentially any computation under any setup of the problem. 

\subsection{The homogeneous setup} \label{subsec:homogeneous}
We now specialize the previous results to the case that  $I_i=I$ and  $J_i=J$ for every $i \in [M]$, where the quantities in \eqref{KL_A}, \eqref{KL_notA}, \eqref{theta_def}, \eqref{system}   simplify as follows:
\begin{align}\label{K_hat_check_homog} 
	\begin{split}
		I_{A} &= I^*_{A} \equiv I,  \;  J_{A} = J^*_{A} \equiv J, \\ 
		\hat{K}_{A} &=
		|A|, \; \check{K}_{A} = |A^c|, \; \theta_{A}  = I/J \equiv  \theta, \quad\\
		c^*_i(A) &=
		\begin{cases}
			x_A,  \;  \text{if}  \; \; i \in A,\\
			y_A,  \;  \text{if}  \; \; i \notin A.
		\end{cases}
	\end{split}
\end{align}

$\bullet$  Suppose first that the  number of anomalous sources is  known in advance,  i.e., $\ell=u$. Then, $x_A$ and $y_A$ do not depend on $A$,  we denote them simply by $x$ and $y$ respectively, and present their values in  Table \ref{tab_2}. 
\begin{table}[h]  
	\begin{center} 
		\begin{tabular}{|c|c|c|c|c|c| c| }
			\hline
			& $ (M-\ell) I \geq  J \ell $& $ (M-\ell)  I\leq  J \ell $  \\
			\hline
			$x$  & $\min \{K/\ell, 1\}$  &  $(K-M+\ell)^+ /\ell$   \\
			\hline  $y$ &	$(K-\ell)^+ /(M-\ell)$  &  $\min \{K/(M-\ell), 1\}$       \\
			\hline
		\end{tabular}
		\caption{$x\equiv x_A$ and $y\equiv y_A$ when $\ell=u$, $I_i=I$, $J_i=J$ for every $i\in [M]$.}
		\label{tab_2}
	\end{center}
\end{table}

From Table \ref{tab_2} and \eqref{Q1} it follows that  $Q_{A}=M$  for every $A  \in \cP_{\ell, u}$.  Then, from   Propositions \ref{propo:prob_ao} and   \ref{propo:prob_ao2}  it follows that if  $R$ is a  probabilistic sampling rule that satisfies  the conditions of Theorem \ref{theo:prob_ao}, then it samples at each time instance each source that is currently  estimated as anomalous (resp. non-anomalous) with probability $x$ (resp. $y$), i.e., 
\begin{align}  \label{system_unique_homog_known}
	c^R_i(A)&=
	\begin{cases}
		x,  \quad  \text{if}  \; \; i \in A\\
		y,  \quad  \text{if}  \; \; i \notin A,
	\end{cases} \forall \; A \in \cP_{\ell, u}.
\end{align}
Moreover, we observe that  the  first-order asymptotic  approximation to the optimal performance is independent of the true subset of anomalous  sources. Specifically, for every $A \in \cP_{\ell, u}$ we have 
\begin{equation}\label{ao_kn0}
	\cJ_{A}(\alpha, \beta, \ell, u, K) \sim  
	\frac{|\log(\alpha \wedge \beta)|}{x \, I  + y \, J}.
\end{equation}


$\bullet$ \quad  When  the  number of anomalous sources is not known a priori, i.e., $\ell < u$, we focus on  the special case that $r=1$ and $\theta =1$, or equivalently,  $I=J$. Then, the values of  $x_A$ and $y_A$,  presented  in  Table   \ref{tab_1}, do not depend on $I$, and  the optimal asymptotic performance under $\Pro_A$ takes the following form:
\begin{equation}\label{ao_lu2}
	\cJ_{A}(\alpha, \beta, \ell, u, K) \sim  \frac{|\log(\alpha)|}{I} \cdot  \, 
	\begin{cases}
		\max\{ (M-\ell)/K, 1\}  , \quad & \text{if} \quad  |A|=\ell,  \\
		M/K, \quad  & \text{if} \quad \ell<|A|<u, \\
		\max\{u/K, 1\} , \quad  & \text{if} \quad |A|=u.
	\end{cases}
\end{equation}
From  \eqref{Q_equiv}  and \eqref{ao_lu2}  we further obtain
\begin{align}  \label{system_unique_homog}
	Q_{A} &=
	\begin{cases}
		M-\ell, \; & \text{if} \quad |A|=\ell, \\
		M, \; & \text{if} \quad \ell <|A|<u, \\
		u, \; & \text{if} \quad |A|=u.
	\end{cases}
\end{align}
Note that,  in this setup, $Q_A=M$ if and only if  one of the following holds:   
$$\ell<|A|<u, \; \;   |A|=\ell=0, \; \; |A|=u=M.$$
Moreover,   from Theorem \ref{theo:AO} it follows that, 
in each of  these three cases,  asymptotic optimality is achieved  by  any sampling rule, not necessarily probabilistic,   that satisfies the sampling constraint and samples all sources  with the same long-run frequency. This is the content of the following corollary. \\

\begin{corollary} \label{coro:ao_unif}
	Let $R$ be a sampling rule such that  
	\begin{itemize}
		\item   the sampling constraint \eqref{samp_const} holds with $T=T^R$,
		\item $\Pro \left(|\pi^R_i(n) -K/M| > \epsilon \right)$ is summable for every   $i \in [M]$ and  $\epsilon>0$.
	\end{itemize}  
	If $\ell<u$, $r=1$, and $I_i=J_j$ for every $i, j \in [M]$, then $R$ is asymptotically optimal under $\Pro_A$ for every $A \in \cP_{\ell, u}$ with  $\ell<|A|<u$. If also  $\ell=0$ and  $u=M$, then $R$ is asymptotically optimal under $\Pro_A$  for every $A \subseteq [M]$.
\end{corollary}

\begin{IEEEproof}
	This follows directly by Theorem \ref{theo:AO}, in view of Table \ref{tab_1}.\\
\end{IEEEproof}

The conditions of  the previous Corollary are satisfied when $R$ is a probabilistic sampling rule with  $c_i^R(A)=K/M$ for every $i \in [M]$ and $A \in \cP_{\ell, u}$, e.g., when  $R$ is a Bernoulli sampling rule that samples  each  source at each time instance with probability $K/M$, independently of the other sources.  Moreover, in the  setup of Corollary \ref{coro:ao_unif}  it is  quite convenient to find and implement a Chernoff rule (Subsection \ref{subsec:chernoff}). Indeed,  when $K$ is an integer,  the conditions of Corollary \ref{coro:ao_unif}  are satisfied when we  take a \textit{simple random sample} of  $K$ sources at each time instance, i.e., when 
$$
q^R(B;D)= 1/ \binom{M}{K} \quad \text{for all}  \; \; D \in \cP_{\ell, u} \; \; \text{and} \;  \;  B \subseteq [M] \; \;  \text{with} \; \;   |B|=K. 
$$
Finally, in Subsection \ref{subsec:tandem} we will introduce a  non-probabilistic   sampling rule that satisfies the conditions of  Corollary \ref{coro:ao_unif}.

\begin{table}[h] 
	\begin{center} 
		\begin{tabular}{|c|c|c|c|c|c| c| }
			\hline
			& $|A|=\ell=0$& $|A|=u=M$ &  $\ell < |A| < u$ &$0<\ell = |A|< u$ & $\ell < |A|=u < M$ \\
			\hline
			$x_A$  & 0 &  $K/M$ & $K/M$ & 0 & 	$\min\{K/u, 1\}$ \\
			\hline  $y_A$ &	$K/M$ & 0 & $K/M$   & $\min\{ K/(M-\ell), 1\}$ &  0 \\
			\hline
		\end{tabular}
		\caption{$x_A$ and $y_A$ when $\ell<u$, $r=1$, $I_i=J_j$ for every $i, j \in [M]$.}
		\label{tab_1}
	\end{center}
\end{table}



\subsection{A heterogeneous example} \label{subsec:heterogeneous}
We end this section by considering a setup where $M$ is an even number, $\ell < M/2 < u$, $r=1$,  and 
\begin{equation} \label{I_hetero}
	I_{i}=  J_i=
	\begin{cases}
		I, \; & 1 \leq i \leq M/2,\\
		I/\phi\; & M/2<  i \leq M,
	\end{cases}
\end{equation}
for some $\phi  \in (0,1]$. Moreover, we focus on the case that the  subset of anomalous sources is of the form $A=\{1, \ldots, |A|\}$. Then,  the optimal asymptotic performance under $\Pro_A$ takes the following form:
\begin{itemize}
	\item when  $|A|=\ell$,
	\begin{equation} \label{ao_lu22}
		\cJ_{A}(\alpha, \beta, \ell, u, K) \sim \frac{|\log \alpha|}{I/\phi} \max\left\{ \frac{(\phi+1) M/2 -\ell}{K} ,1 \right\},
	\end{equation} 
	\item when $\ell<|A|<u$,
	\begin{equation} \label{ao_lu33}
		\cJ_{A}(\alpha, \beta, \ell, u, K) \sim \frac{|\log \alpha|}{I} \, \max\left\{ \frac{(\phi+1) M/2}{K} , 1 \right\},
	\end{equation} 
	\item when $|A|=u$, 
	\begin{equation}\label{ao_lu44}
		\cJ_{A}(\alpha, \beta, \ell, u, K) \sim  \frac{|\log \alpha|}{I} 
		\max\left\{ \frac{(1-\phi) (M/2)+ \phi u}{K},  1 \right\} . 
	\end{equation} 
\end{itemize}
From these expressions and \eqref{Q_equiv} we  also obtain 
\begin{equation}
	Q_{A}=
	\begin{cases}
		(\phi+1) M/2 -\ell, \quad & |A|=\ell, \\
		(\phi+1) M/2, \quad & \ell <|A|<u, \\
		(1-\phi) (M/2)+ \phi u, \quad  &|A|=u,
	\end{cases}
\end{equation}
and we note that  $Q_A$ is  always strictly smaller than $M$  in this setup as long as 
$\phi<1$.

\section{Non-probabilstic sampling rules} \label{sec:other}
In this section, we discuss certain non-probabilistic sampling rules. 

\subsection{Sampling in tandem}  \label{subsec:tandem}
Suppose that $K$ is an integer and consider the straightforward sampling approach according to which the sources are sampled in tandem, $K$ of them at a time. Specifically,   sources $1$ to $K$ are sampled at time $n=1$, 
and if $2K \leq M$, then sources $K+1$ to $2K$ are sampled at time $n=2$, whereas if $2K>M$, then sources $K+1$ to $M$ and $1$ to  $2K-M$ are sampled at time $n=2$,   etc. In this way,  each source is sampled exactly $K$ times in an interval  of the form $((m-1)M, m M]$, where $m \in \bN$, to which we will refer as a \textit{cycle}.   This sampling rule satisfies the conditions of Corollary \ref{coro:ao_unif}, which means that 
in the special case that $\ell<u$, $r=1$, and $I_i=J_j$  for every $i, j \in [M]$, it achieves asymptotic optimality under $\Pro_A$ when $\ell<|A|<u$, and for every $A \subseteq [M]$ when $\ell=0$ and $u=M$. 

In general, by the formula for the optimal asymptotic performance under  full sampling, which is obtained by the lower bound of Theorem \ref{theo:lowerbound} with $K=M$, we can see that if sampling is terminated  at a time that is a multiple of $M$,  the expected \textit{number of  cycles} until stopping is, to a first-order asymptotic approximation, equal to $\cJ(\alpha, \beta, \ell,  u, M) /K$. Since each cycle is  of length $M$, the expected \textit{time} until stopping is, again to a first-order asymptotic approximation, equal to  $ (M/K) \cJ(\alpha, \beta, \ell,  u, M)$. Thus,  the  \textit{asymptotic relative efficiency} of this sampling  approach  can be defined as follows:
\begin{equation} \label{ARE}
	{\sf ARE}:=\frac{M}{K} \lim_{\alpha, \beta \to 0}  \frac{\cJ_A(\alpha, \beta, \ell, u, M)}{\cJ_A(\alpha, \beta,\ell, u, K)},
\end{equation}
where the limit is taken so that \eqref{r} holds when $\ell<u$. 

Consider in particular the homogeneous setup of Subsection \ref{subsec:homogeneous}, where \eqref{K_hat_check_homog} holds. When
$\ell=u$,   by \eqref{system_unique_homog_known} and \eqref{ao_kn0} we have:
\begin{itemize}
	\item if $\theta \geq  \ell/ (M-\ell)$,
	\begin{equation}\label{ARE_known1}
		{\sf ARE}= \frac{\theta}{1+\theta} \, \frac{M}{\max\{\ell,K\} } + \frac{1}{1+\theta} \frac{(1-\ell/K)^+}{1-\ell/M},
	\end{equation} 
	\item 
	if $\theta <  \ell/ (M-\ell)$,
	\begin{equation}\label{ARE_known2}
		{\sf ARE}= \frac{1}{1+\theta} \, \frac{M}{\max\{M-\ell,K\} } + \frac{\theta}{1+\theta} \frac{(1-(M-\ell)/K)^+}{\ell/M}.
	\end{equation} 
\end{itemize}
When $\ell<u$ and $r=\theta=1$,   by \eqref{ao_lu2} we have 
\begin{equation} \label{ARE_unknown}
	{\sf ARE}= \begin{cases}
		M / \max\{ M-\ell, K \} , \quad & \text{if} \quad  |A|=\ell,  \\
		1, \quad  & \text{if} \quad \ell<|A|<u ,\\
		M / \max\{u, K \}  \quad  & \text{if} \quad |A|=u.
	\end{cases}
\end{equation}
On the other hand,  in the heterogeneous setup  of Subsection \ref{subsec:heterogeneous}, by 
\eqref{ao_lu22}, \eqref{ao_lu33}, \eqref{ao_lu44} we have:
\begin{equation} \label{ARE_hetero}
	{\sf ARE}=\begin{cases}
		M/ \max\{ (\phi+1) M/2-\ell, \; K\}, \; & \text{if} \quad  |A|=\ell,  \\
		M/ \max\{ (\phi+1)M/2, \; K\}, \;  & \text{if} \quad \ell<|A|<u ,\\
		M/ \max\{ (1-\phi)M/2 +\phi u, \; K\} \;  & \text{if} \quad |A|=u.
	\end{cases}
\end{equation}

\subsection{Equalizing empirical and limiting sampling frequencies}
We next  consider a   sampling  approach, which  has been applied to a general controlled sensing problem  \cite{Aditya_2021}, as well as to a bandit problem \cite{garivier2016optimal}, and we show that  not only  it may not achieve asymptotic optimality in the sequential anomaly detection problem, but  it may  even lead to a detection  procedure that \textit{fails to terminate with positive probability}.   

To be more specific, we consider the  homogeneous setup of Subsection \ref{subsec:homogeneous} with $K=1$. In this setup,    a probabilistic sampling rule  that satisfies the conditions of Theorem  \ref{theo:prob_ao}  samples a source in 
$D$ (resp.  $D^c$) with probability $x_{D}$ (resp. $y_{D}$), whenever $D \in \cP_{\ell, u}$ is the subset of sources currently identified as  anomalous.  The sampling rule $R$ that  we consider  in this Subsection is not probabilistic, as it  requires knowledge  of not only  the currently  estimated subset of anomalous sources, but also of the current  empirical sampling  frequencies. Specifically,  if $D$ is the subset of sources currently estimated as anomalous,  for every source in $D$  (resp. $D^c$)  it computes the  distance between its  current empirical sampling frequency  and  $x_{D}$ (resp. $y_{D}$), and it samples next a source for which this distance is the maximum. That is,  for every $n \in \bN$ and $D \in \cP_{\ell, u}$ we have  
\begin{align}\label{eqem}
	\begin{split}
		R(n+1) &\in \text{argmax} \left\{ |\pi_i^R(n)- x_D|, |\pi_j^R(n)- y_D| : 
		i \in D,  \; j \notin D \right\} \\
		&\text{on} \quad  \{\Delta^R_n=D\}.
	\end{split}
\end{align}
Without loss of generality, we also assume that each source  has  positive probability to be sampled  at the first time instance, i.e., 
\begin{align}\label{eqem2}
	\Pro_{A}(i \in R(1) )>0 \quad \forall  \; i \in [M], \; A \in \cP_{\ell, u}.
\end{align}

\begin{proposition}
	Consider the homogeneous setup of Subsection \ref{subsec:homogeneous}  with $K=1$ and let $R$ be sampling rule that satisfies \eqref{eqem}-\eqref{eqem2}.  Suppose further that there is only one anomalous source, i.e.,
	that  the subset of anomalous source, $A$,   is a singleton, and also that  $x_A+y_A<1$.
	Then, there is  an event  of positive probability under $\Pro_A$ on which 
	\begin{enumerate}
		\item[(i)]   the same source is sampled at every time instance,
		\item[(ii)]  and if also $\ell=0$ and $u=M$,  $T^R$  fails to terminate for any selection of its thresholds. 
	\end{enumerate}
\end{proposition}

\begin{IEEEproof}
	If (i) holds,  there is  an event  of positive probability under $\Pro_A$ on which   all  LLRs but one are always equal to 0. Thus,  (ii)  follows directly from (i) and  the fact that when 
	$\ell=0$ and  $u=M$,   the stopping rule \eqref{intersection} requires that all LLRs be non-zero upon stopping.  Therefore, it remains to prove (i). Without loss of generality, we set  $A=\{1\}$.  The event $\{ R_1(1)=1\}$ has positive probability under $\Pro_A$, by \eqref{eqem2},  and so does the event  \begin{equation} \label{event0}
		\Gamma \equiv \left\{\sum_{m=1}^n g_1(X_1(m)) >0 \; \; \; \forall  \; n \in \bN \right\},
	\end{equation}
	since  $\{g_1(X_1(n)): n \in \bN\}$ is an iid sequence  with expectation  $I>0$ under $\Pro_A$ (see, e.g., \cite[Proposition 7.2.1]{resnick1992adventures}). Since these two events are also independent, their  intersection  has  positive probability under $\Pro_A$. Therefore, it suffices to show that only source 1 is sampled on $\Gamma \cap \{ R_1(1)=1\}$.   Indeed, on this event  sampling starts from source 1, and as a result  the vector of empirical frequencies, $(\pi_1^R(n), \ldots, \pi_M^R(n))$ at time $n=1$ is $(1,0, \ldots, 0)$. Moreover,  the estimated  subset of anomalous sources at $n=1$  is  $A=\{1\}$, independently of  whether $\ell<u$ or $\ell=u$.  Since  $x_A+y_A<1$,  or equivalently, 
	$$|\pi_1(1)-x_A|=|1-x_A|= 1-x_A> y_A=  |y_A-0|=|\pi_i(1)-y_A|$$
	for every $i \neq 1$,  by  \eqref{eqem} it follows that  source $1$ is sampled again at time $n=2$, and as a result  the vector of empirical frequencies  at time $n=2$ remains $(1,0,\ldots,0)$. 
	Applying the same reasoning as before, we conclude  that source 1 is sampled again at  time $n=3$. The same argument can be repeated indefinitely, and this  proves (i).   \\
\end{IEEEproof}

\noindent \underline{\textbf{Remark:}} When
$\ell=0$, $u=M$,   each  $f_{0i}$  is exponential with rate $1$, and each $f_{1i}$ is exponential with  rate $\lambda >0$,  the conditions of the previous proposition are satisfied as long as $M>2$. Indeed, in this case we have 
\begin{align*}
	\theta&=I/J=\frac{-\log(\lambda)+\lambda-1}{\log(\lambda)+1/\lambda -1}, \\
	x_A &=\frac{1}{1+(M-1)\theta}, \quad y_A= \theta x_A  
\end{align*}
and as a result  $x_A+y_A<1 \Leftrightarrow M>2$.


\subsection{Sampling based on the ordering of the LLRs} \label{sub:deterministic}

A different,  non-probabilistic sampling approach, which goes back to \cite[ Remark 5]{chernoff1959}, suggests  sampling at each time  instance the sources with the smallest, in absolute value,  LLRs among those estimated as anomalous/non-anomalous.  Such a sampling rule was proposed  in \cite{Cohen2015active}, in the  homogeneous setup of Subsection \ref{subsec:homogeneous}, under the assumption that the  number of anomalous sources is known a priori, i.e., $\ell=u$. An extension of this rule in the heterogeneous setup was studied  in  \cite{huang2017active}, under the assumption that $\ell=u$ and $K=1$. A similar sampling rule,  that also has  a randomization feature, was  proposed in   \cite{Tsopela_2019}  when  $\ell =u$, as well as in the case of no prior information, i.e., when   $\ell=0,u=M$.   For each of them,  the criterion of Theorem \ref{theo:AO} can be applied  to establish their asymptotic optimality. Its  verification,  however, is  a quite  difficult task that we plan to  consider in the future.

\section{Simulation study}\label{sec:simulations}
In this section we present the results of a simulation study where we  compare  two   probabilistic sampling rules, Bernoulli (Section \ref{sec:probabilistic})  and Chernoff (Subsection \ref{subsec:chernoff}),  between them and against   sampling in tandem (Subsection \ref{subsec:tandem}).  Throughout this section,  for every $i \in [M]$ we have  $f_{0i}=\mathcal{N}(0,1)$ and $f_{1i}=\mathcal{N}(\mu_i,1)$,  i.e., all observations from source $i$ are Gaussian with  variance $1$   and mean equal to  $\mu_i$ if the source is anomalous and 0 otherwise,   and as a result  $I_{i}=J_{i}=(\mu_{i})^2 /2$.  We consider a  homogeneous setup where 
\begin{equation}\label{mu_homog}
	\mu_{i}=\mu, \quad  \forall \;  i \in [M],
\end{equation}
as well as  a  heterogeneous setup  where
\begin{equation}\label{mu_hetero}
	\mu_{i}=
	\begin{cases}
		\mu,\quad & 1 \leq i\leq M/2\\
		2 \mu ,\quad & M/2 < i \leq M,
	\end{cases}
\end{equation}
in which case  \eqref{I_hetero} holds with $I=\mu^2/2$ and $\phi=0.25$. In both setups we set $\alpha=\beta$, $M=10$,  $K=5$, $\ell =1$,  $u=6$, $\mu=0.5$, we assume that  the  subset  of anomalous sources is  of the form $A=\{1, \ldots, |A|\}$, and  consider different values for its size. 
The two probabilistic sampling rules are designed so that   \eqref{system_epsilon}  holds for every $A \in \cP_{\ell, u}$. As a result,  by Corollary \ref{coro:prob_ao} it follows that in both setups they are asymptotically optimal under $\Pro_A$ for every $A \in \cP_{\ell, u}$.  On the other hand, by Corollary \ref{coro:ao_unif} it follows that sampling in tandem is asymptotically optimal under $\Pro_A$ only in the  homogeneous setup \eqref{mu_hetero} and when  $l<|A|<u$ (since $0<l<u<M$).  

In Figure \ref{fig:nonasym} we  plot the  expected value of the  stopping time that is induced by  each  of the three sampling rules against the  true number of anomalous sources.  Specifically,  in Figure \ref{fig:homog_nonasym} we consider the homogeneous setup \eqref{mu_homog} and in Figure \ref{fig:nonhomog_nonasym} the heterogeneous setup \eqref{mu_hetero}. In all cases, the  thresholds  are chosen, via Monte Carlo simulation, so that the familywise  error  probability of each kind is   (approximately) equal to $\alpha=\beta=10^{-3}$.  The  Monte Carlo standard error that corresponds to each estimated expected value is approximately $10^{-2}$. From Figure \ref{fig:nonasym}  we can see that the performance implied by the  two probabilistic sampling rules is always essentially  the same. Sampling in tandem performs significantly worse in all cases apart from the  homogeneous setup with $\ell<|A|<u$, where all three sampling rules lead to  essentially  the same performance. We note also that, in both setups and for all three sampling rules,  the expected time until stopping  is  much smaller when  the number of anomalous sources is equal to either  $\ell$ or $u$ than when it is between  $\ell$ and $u$. 

\begin{center}
	\begin{figure}[h]
		\subfloat[Homogeneous case]{
			\includegraphics[width=0.5\linewidth]{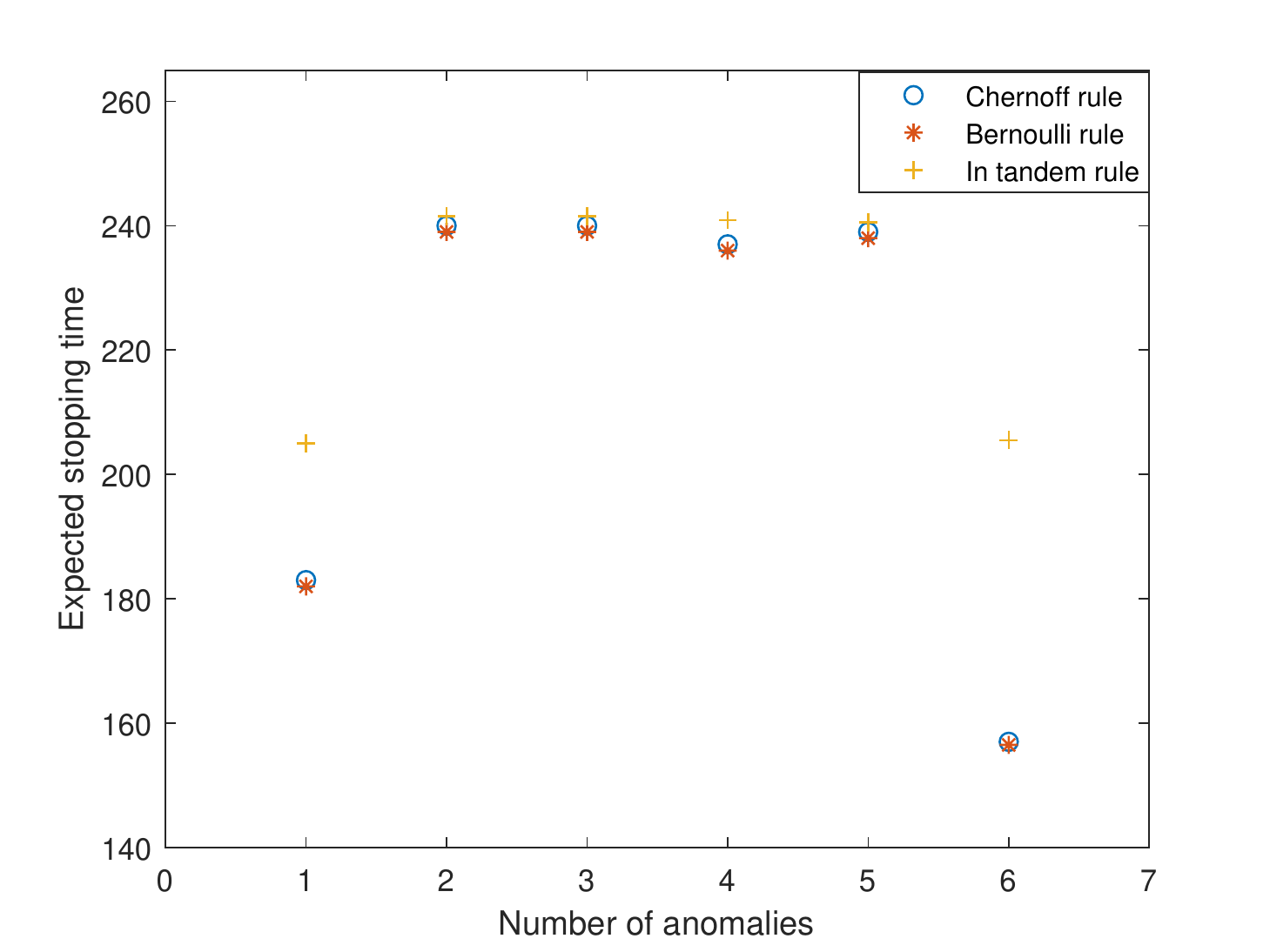}  
			\label{fig:homog_nonasym}
		}
		\subfloat[Heterogeneous case]{
			\includegraphics[width=0.5\linewidth]{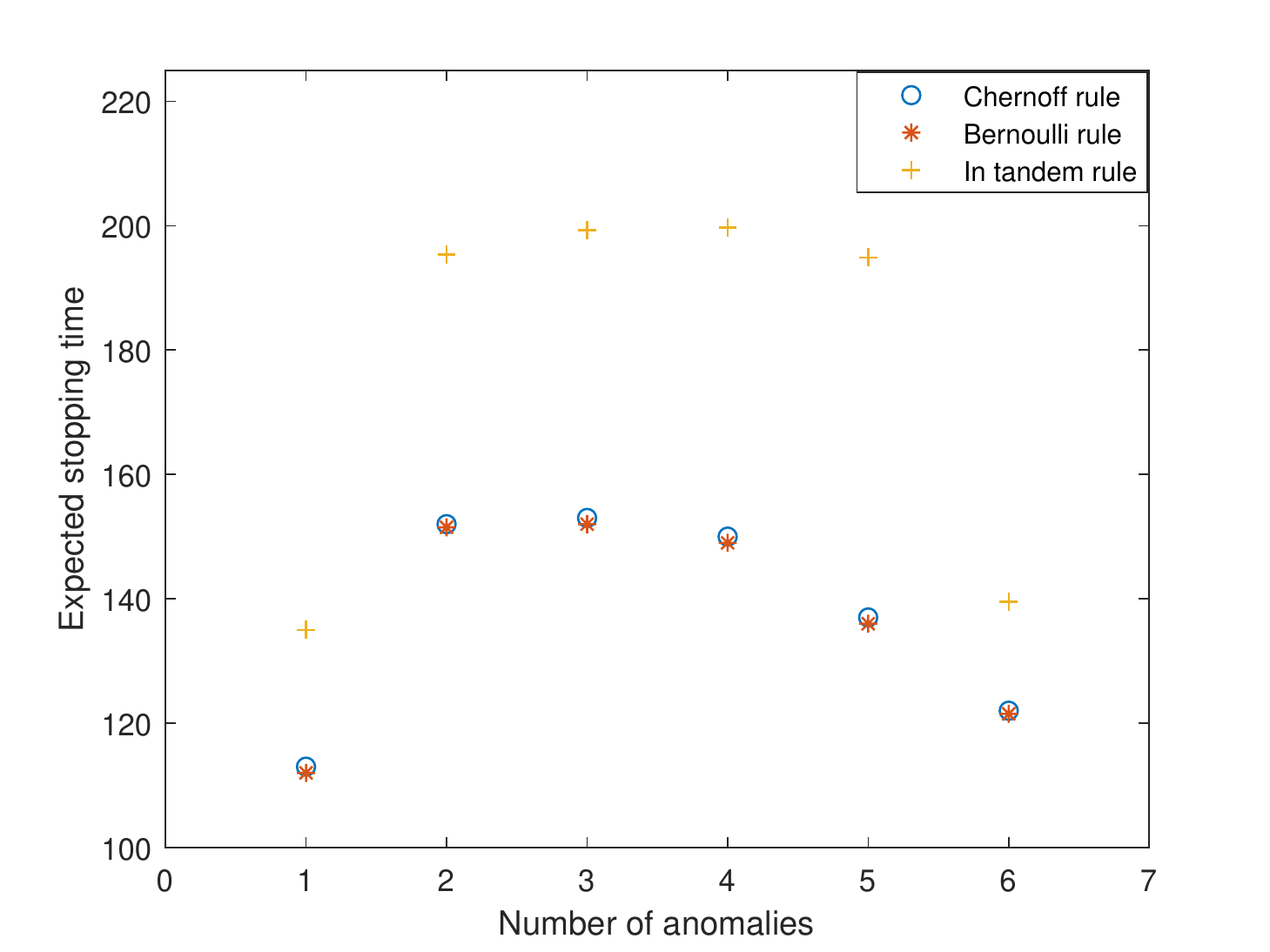}
			\label{fig:nonhomog_nonasym} 
		} 
		\captionsetup{justification=raggedright, singlelinecheck=false}
		\caption{Expected value of the stopping time that corresponds to each of the three sampling rules versus the number of anomalous sources ($\alpha=\beta=10^{-3}$).}
		\label{fig:nonasym} 
	\end{figure}
\end{center}

\begin{center}
	
	\begin{figure}[h]
		\subfloat[$|A|=1$]{
			\includegraphics[width=0.51\linewidth]{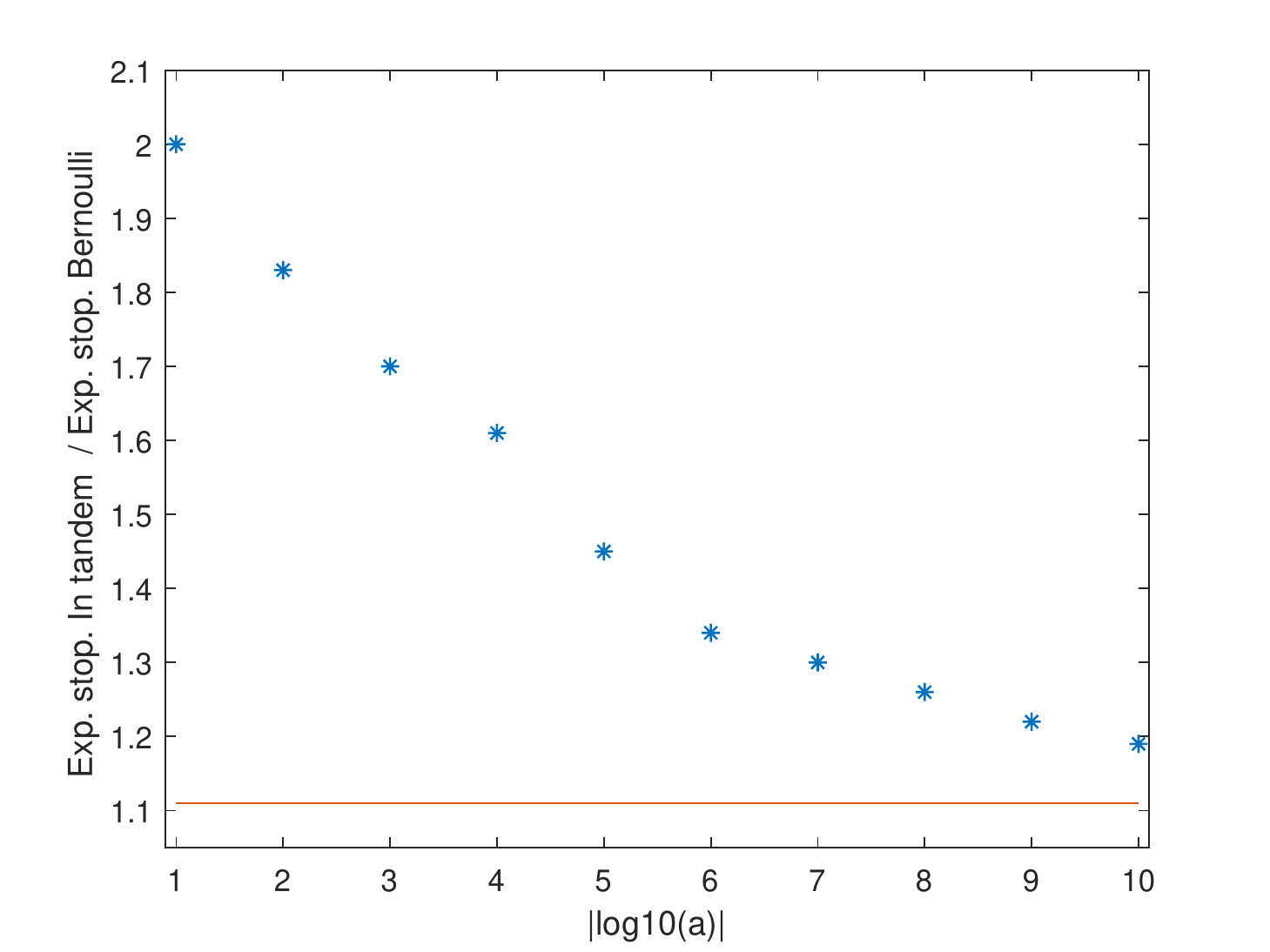}  
			\label{fig:m0i}
		}
		\subfloat[$|A|=6$]{
			\includegraphics[width=0.51\linewidth]{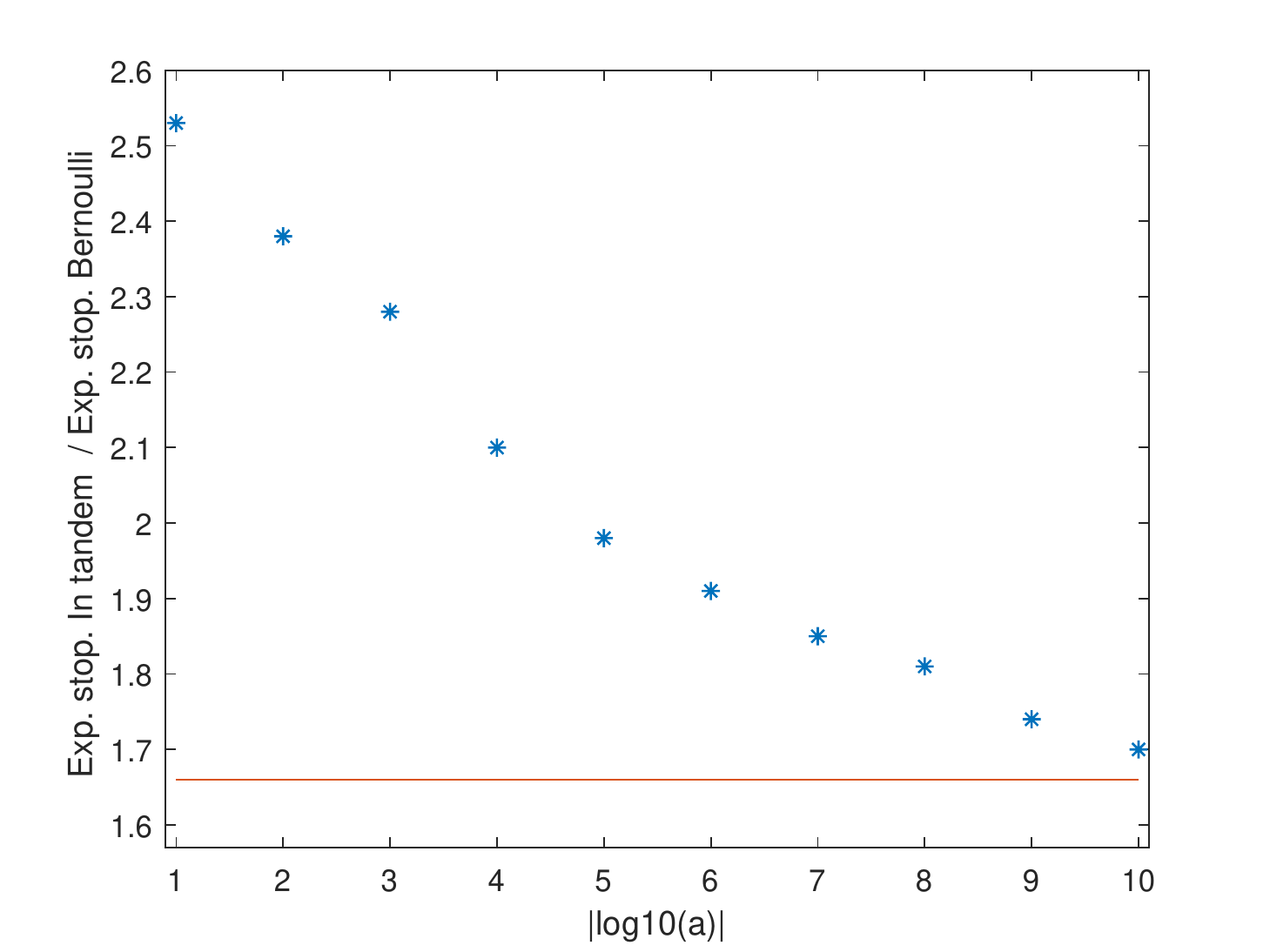}  
			\label{fig:m_n0i}
		} 
		\caption{In both graphs the $x$-axis  represents $|\log_{10}(\alpha)|$ and the $y$-axis the ratio of the expected value of the stopping time implied by sampling in tandem over that  implied by the  asymptotically optimal Bernoulli rule in the  homogeneous setup \eqref{mu_homog}.  The horizontal  line refers to the value of \eqref{ARE_unknown}.}
		\label{fig:homog_asymi}
	\end{figure}
\end{center}

In Figures \ref{fig:homog_asymi} and \ref{fig:nonhomog_asymi}
we  plot the ratio of the expected value of the stopping time induced by sampling in tandem over that induced by the Bernoulli sampling rule  against $|\log_{10}(\alpha)|$  as  $\alpha$ ranges from $10^{-1}$ to $10^{-10}$. (We do not present the corresponding results for the  Chernoff rule, as they are almost identical).  Specifically, in Figure \ref{fig:homog_asymi} we consider the  homogeneous setup \eqref{mu_homog}  when the number of anomalous sources is $1$ and $6$, whereas in Figure \ref{fig:nonhomog_asymi} we consider the   heterogeneous setup \eqref{mu_hetero} when the  number of anomalous sources is $1, \, 3,$ and  6.  For each value of $\alpha$, the thresholds are selected  according to  \eqref{thresholds_gi} and each expectation is computed using $10^4$ simulation runs. The standard error for each estimated expectation  is approximately  $1$, whereas  the   standard error for each ratio is  approximately  $10^{-2}$ in  the homogeneous setup and $3 \cdot 10^{-2}$ in the heterogeneous setup.  Moreover, in each case we plot the limiting value of this ratio, which is the limit defined in \eqref{ARE}.  In the homogeneous case this is given by  \eqref{ARE_unknown}  and  is  equal to 
\begin{equation} 
	\begin{cases}
		10/ \max\{ 10-1 ,\; 5\} =10/ 9 \approx 1.1 
		& \;  \text{when} \quad  |A|=1,  \\
		10/ \max\{ 6,  \; 5\} =10/ 6 \approx 1.6 \; &
		\;  \text{when} \quad  |A|=6.  \\
	\end{cases}
\end{equation}
In the heterogeneous case   it is given by \eqref{ARE_hetero} and  is equal to 
\begin{equation} 
	\begin{cases}
		10/ \max\{ (1+0.25) \cdot 5-1, \; 5\} =10/ 5.25 \approx 1.9 \; & \text{when} \quad  |A|=1,  \\
		10/ \max\{ (1+0.25) \cdot 5, \; 5\}= 10/ 6.25 \approx 1.6  \;  & \text{when} \quad 1<|A|<6 ,\\
		10/ \max\{ (1-0.25) \cdot 5 + 0.25 \cdot 6, \; 5\} =10/ 5.25 \approx 1.9 \;  & \text{when} \quad |A|=6.
	\end{cases}
\end{equation}
From Figures \ref{fig:homog_asymi} and \ref{fig:nonhomog_asymi} we can see that, in all cases, the efficiency loss  due to  sampling in tandem is (much)  larger than the one suggested by the corresponding asymptotic relative efficiency when  $\alpha$ is not (very) small. 

\begin{center}	
	\begin{figure}[h]
		\subfloat[$|A| =1$]{
			\includegraphics[width=0.33\linewidth]{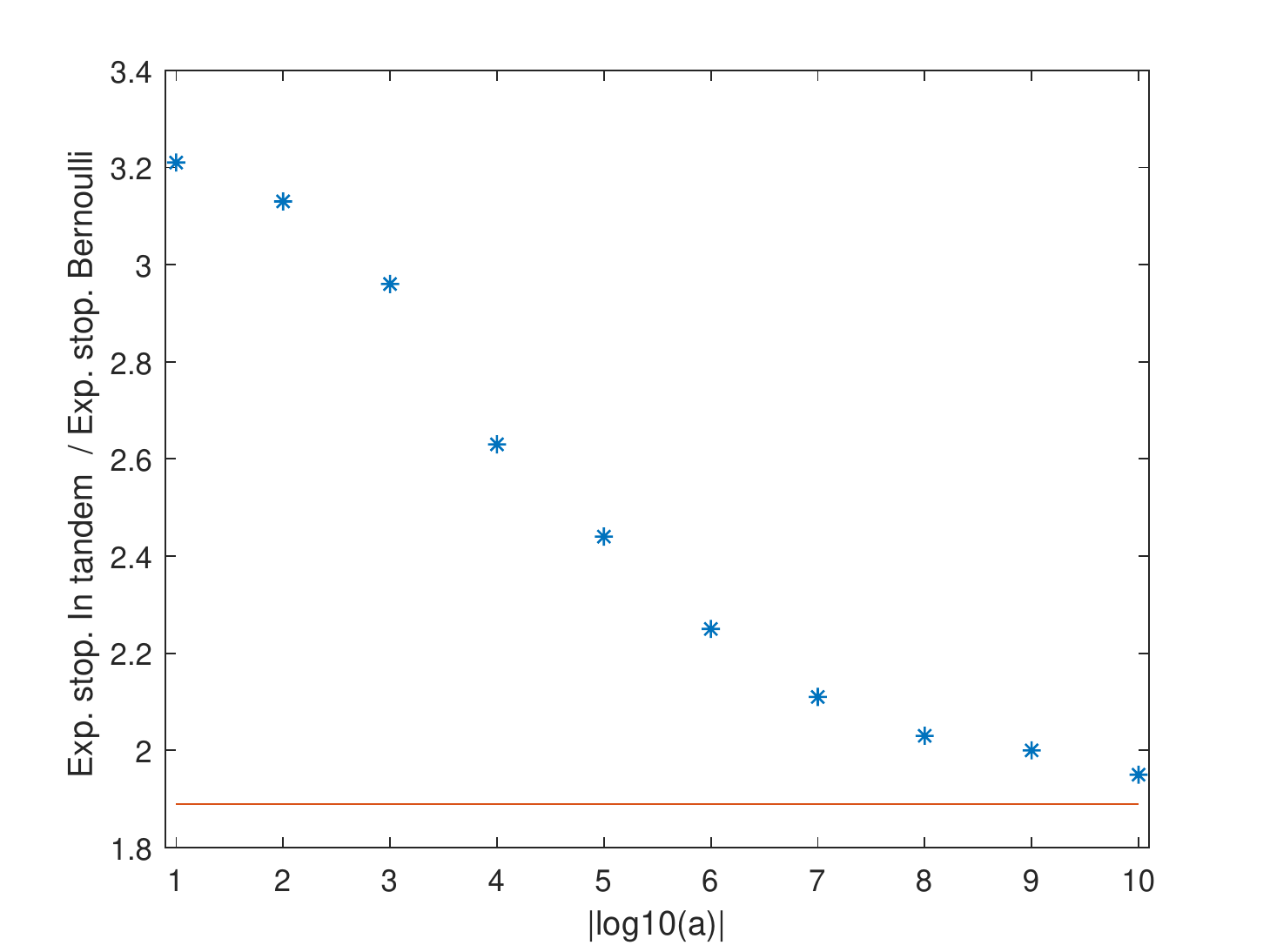}  
			\label{fig:mi}
		}
		\subfloat[$|A| = 3$]{
			\includegraphics[width=0.33\linewidth]{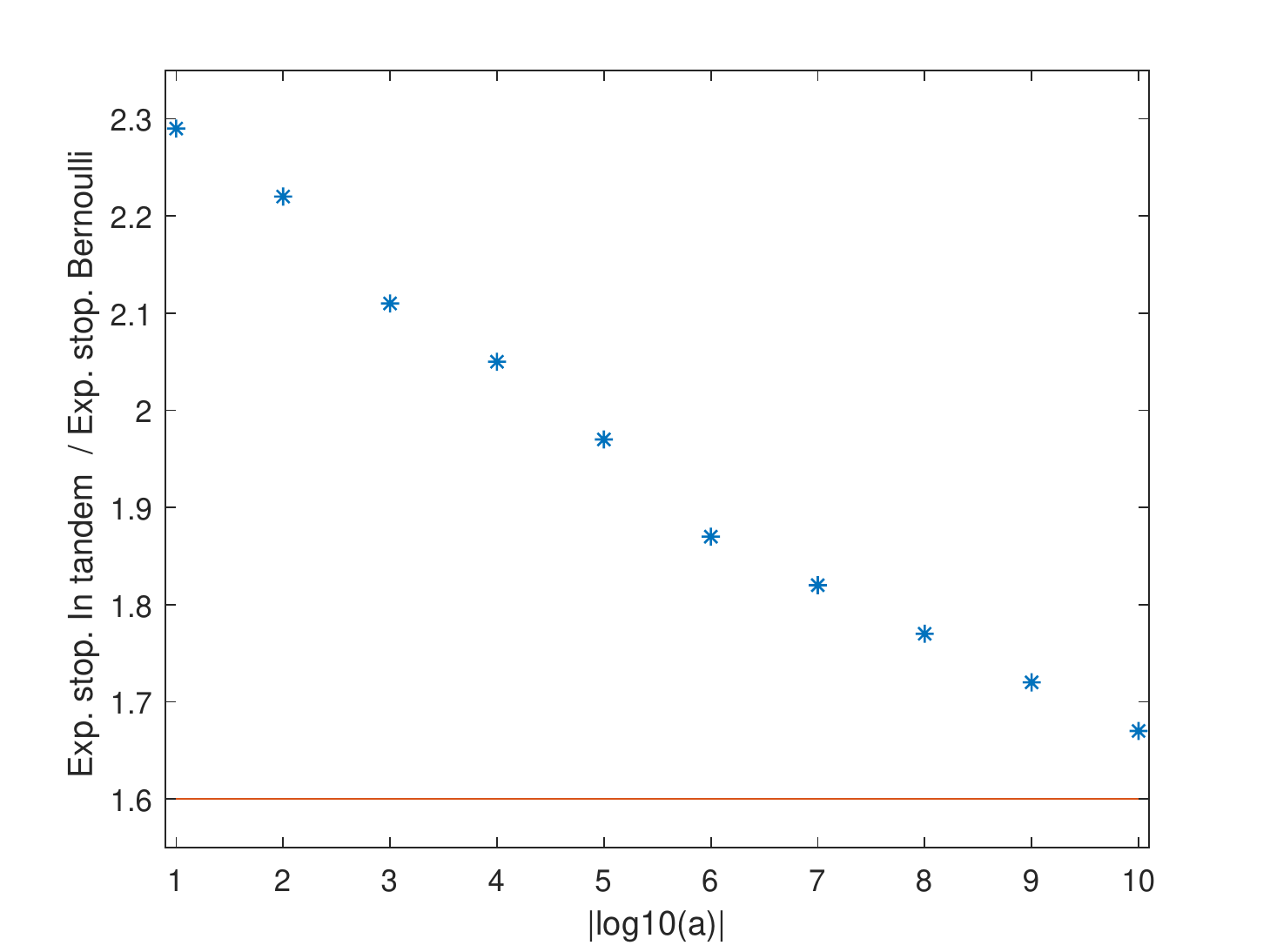}  
			\label{fig:luii}
		}
		\subfloat[$|A|= 6$]{
			\includegraphics[width=0.33\linewidth]{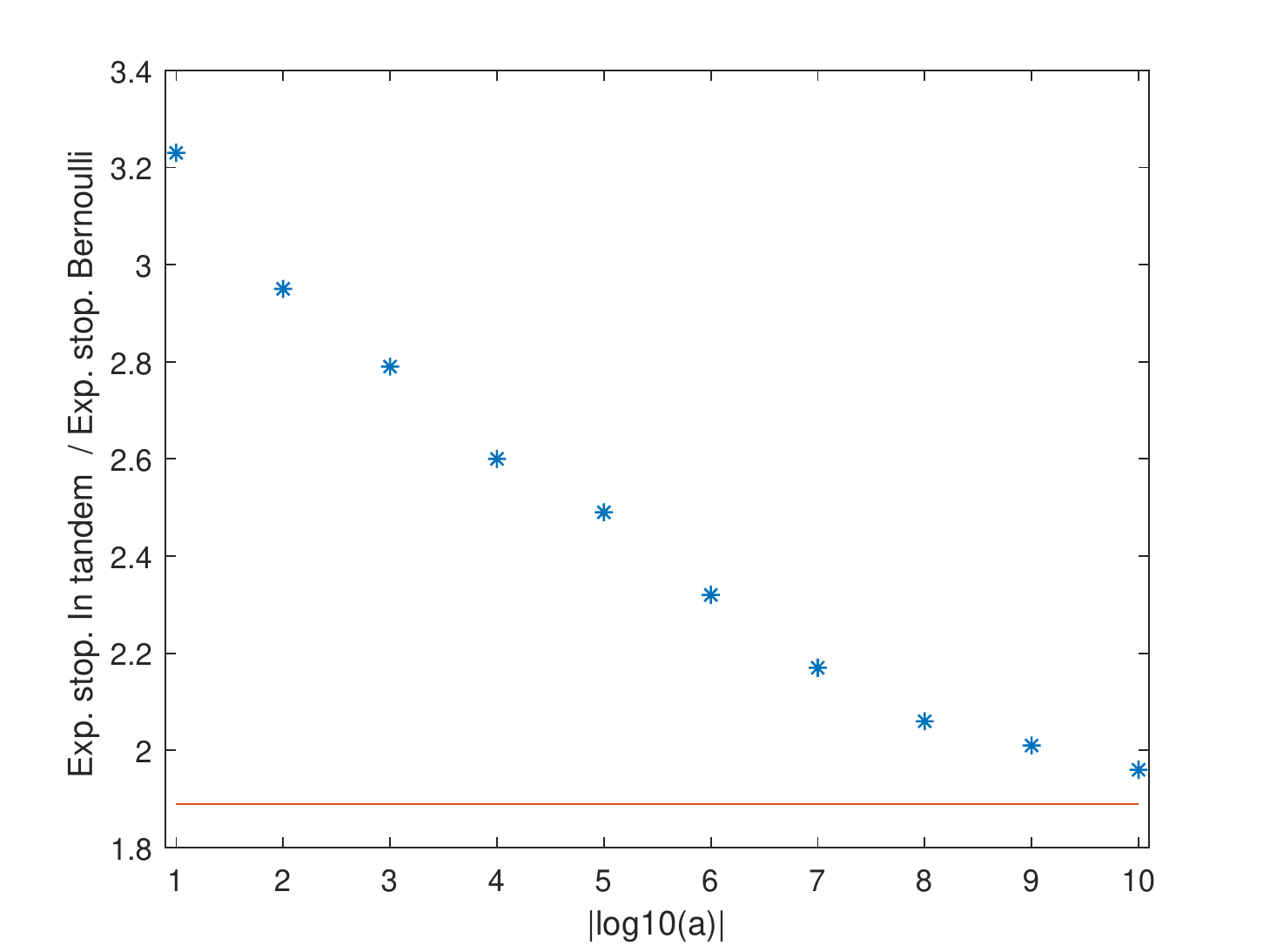}  
			\label{fig:m_ni}
		} 
		\caption{In all graphs the $x$-axis  represents $|\log_{10}(\alpha)|$ and the $y$-axis the ratio of the expected value of the stopping time implied by sampling in tandem over  that implied by the  asymptotically optimal Bernoulli rule in the  heterogeneous setup \eqref{mu_hetero}.  The horizontal  line refers to the value of \eqref{ARE_hetero}.}
		\label{fig:nonhomog_asymi} 
	\end{figure}
\end{center}

\section{Conclusions and extensions}  \label{sec:conclusions}
In this paper we propose a novel formulation of the sequential anomaly detection problem with sampling constraints, in which arbitrary, user-specified  bounds are  assumed on the number of anomalous sources,  the probabilities of at least one false alarm and at least one missed detection are controlled below distinct tolerance levels,  and the  number  of sampled sources per time instance is not necessarily fixed.  We obtain  a general criterion for achieving the optimal expected time until stopping,  to a first-order asymptotic approximation as the error probabilities go to 0,  as long as the log-likelihood ratio statistic of each observation has a   finite \textit{first} moment. We show that asymptotic optimality is achieved, simultaneously under every possible subset of anomalous sources,  for  any version of the proposed problem, using the \textit{unmodified} sampling rule in \cite{Bessler1960_I}, to which we refer  as \textit{Chernoff}, but also using a novel sampling rule whose implementation  requires minimal computations,  to which we refer  in this work as \textit{Bernoulli}. Despite their very different computational requirements, these two rules are found in simulation studies to lead to essentially the same performance.
On the other hand, 
it has been  shown in  simulation studies  under various setups \cite{Cohen2015active,huang2017active,Tsopela_2019},   that the Chernoff rule leads to   significantly worse performance  in practice  compared to non-probabilistic sampling rules  as the ones   discussed in Subsection  \ref{sub:deterministic}.   The study of  sampling rules of this nature in the general sequential anomaly detection framework we propose in this work is an interesting problem that we plan to consider in the future. Another interesting problem is that of  establishing stronger than first-order  asymptotic optimality, which has been solved  in certain special cases of the  sequential multi-hypothesis testing problem with controlled sensing \cite{Keener_1984,Lalley_Lorden_1986}.  

The results of this paper  can  be shown, 
using the techniques in \cite{Bartroff_and_Song_2020},  to remain valid for a variety of error metrics beyond the  familywise error rates that we consider in this work, such as the false discovery rate  and the false non-discovery rate.  However, this is not the case  for  the generalized familywise error rates proposed in \cite{hommel1988controlled,LehmannRomano2005}, for which different policies and a different analysis is required. These error metrics have been considered in \cite{De_and_Baron_2012b,Bartroff_and_Song_2014,Song_and_Fellouris_2019} in the case of full sampling, whereas  certain results in the presence of sampling constraints have been presented in \cite{Tsopela_2020}.

The results in this work can also be generalized in a natural  way when the sampling cost varies per source, as in \cite{Cohen2019nonlinearcost}, or when  the two hypotheses in  each source are not completely specified, as it is done for example in   \cite{Song_and_Fellouris_2019} in the case of full sampling.  Another potential generalization is the removal of the assumption that the acquired observations are  conditionally independent of the past given the current sampling choice, as it is done in  \cite{nitinawarat_controlled_2015}  in a general controlled sensing setup.  Finally, another direction of interest is a setup where the focus is on  the dependence structure of the sources rather than their marginal distributions, as for example  in  \cite{heydari2016quickest}.

\section{Acknowledgments} 
The authors would like to thank Prof. Venugopal Veeravalli  for stimulating discussions.

\appendices
\renewcommand{\thetheorem}{ \Alph{section}.\arabic{theorem}}
\renewcommand{\thelemma}{ \Alph{section}.\arabic{lemma}}

\section{} \label{appen:A}

In this Appendix we state and prove two auxiliary lemmas that are used in the proofs of various results of this paper. Specifically, we fix $A \in \cP_{\ell, u}$,  and   for any  $i \in [M]$,  $n \in \bN$ and any sampling rule  $R$ we set  
\begin{align} \label{tilde}
	\begin{split}
		\widetilde{\Lambda}^R_i(n) &:= \widetilde{\Lambda}^R_i(n-1) + \Bigl( g_{i}(X_i(n))   - \Exp_A[ g_i(X_i(n))] \Bigr) \, R_i(n),\\
		\widetilde{\Lambda}^R_i(0) &:= 0,
	\end{split}
\end{align} 
and  comparing with \eqref{LLR} we observe that 
\begin{align} \label{decompose}
	\widetilde{\Lambda}^R_i(n) &= 
	\begin{cases}
		\Lambda^R_i(n) - I_i \, N^R_i(n), \quad i \in A \\ 
		\Lambda^R_i(n) + J_i \, N^R_i(n), \quad i \notin A. 
	\end{cases}
\end{align}

\begin{lemma}\label{lem:A1}
	Let  $R$ be an arbitrary sampling rule,  $i \in A$,  $j \notin A$,  $\rho \in (0,1]$, and  $\epsilon>0$. Then, the sequences 
	\begin{align*}
		&\Pro_{A}\left( \frac{1}{n} \Lambda^R_i(n)< 
		\rho I_i-\epsilon, \; \pi_i^R(n) \geq \rho \right), \\
		&\Pro_{A}\left( \frac{1}{n}  \Lambda^R_j(n)>-\rho J_j+\epsilon, \; \pi_j^R(n) \geq \rho \right), \\
		&\Pro_{A}\left(\frac{1}{n} \Lambda^R_{ij}(n)< \rho I_i-\epsilon, \; \pi^R_i(n) \geq \rho \right) \\
		&\Pro_{A}\left( \frac{1}{n}  \Lambda^R_{ij}(n)< \rho J_j-\epsilon, \; \pi^R_j(n) \geq \rho \right)
	\end{align*} 
	are exponentially decaying. 
\end{lemma}

\begin{IEEEproof}
	We  prove the result for the third probability, as the proofs  for the other ones are similar. To lighten the notation, we suppress the dependence on $R$ and  we  write 
	$\widetilde{\Lambda}_{ij}(n), \widetilde{\Lambda}_i(n),  \pi_i(n), \cF_n$ instead of  
	$\widetilde{\Lambda}^R_{ij}(n), \widetilde{\Lambda}^R_i(n), \pi_i^R(n), \cF^R_n$. By \eqref{decompose}, for every $n \in \bN$  we have 
	\begin{align*}
		\Lambda_{ij}(n) &= 
		\Lambda_{i}(n) -\Lambda_{j}(n) \\
		&= 
		\widetilde{\Lambda}_{i}(n) -\widetilde{\Lambda}_{j}(n) + n (I_i \pi_i(n)+ J_j \pi_j(n)),
	\end{align*}
	which  shows that if $\pi_i(n) \geq \rho$, then 
	$
	\Lambda_{ij}(n)  \geq  
	\widetilde{\Lambda}_{i}(n) -\widetilde{\Lambda}_{j}(n) + n \rho I_i. $
	Thus, for every $n \in \bN$ we have 
	\begin{align*}
		& \Pro_{A}\left( \frac{1}{n}  \Lambda_{ij}(n) < \rho I_i-\epsilon, \; \pi_i(n) \geq \rho \right) \\
		&\leq	 \Pro_{A}\left( \widetilde{\Lambda}_i(n)  - \widetilde{\Lambda}_j(n) <-n \, \epsilon  \right)\\
		&\leq	 \Pro_{A}\left( \widetilde{\Lambda}_i(n) <-n\,  \epsilon/2 \right) +  \Pro_{A}\left(-\widetilde{\Lambda}_j(n) < -n \, \epsilon/2 \right), 
	\end{align*}
	and it suffices to show that the two terms in the  upper bound are exponentially decaying. We  show this only for the first one, as the proof for the second is similar.  For this, we fix $\delta \in (0, \epsilon/2)$ and we observe that, for any  $t>0$ and $n \in \bN$, by Markov's inequality we have
	\begin{align} \label{markov0}
		\begin{split}
			&  \Pro_{A}  \left(  \widetilde{\Lambda}_i(n)<-n\,  \epsilon/2 \right) \\
			&\leq  \exp\{-n\,t\,  \epsilon/2 \}   \, \Exp_A \left[    \exp \left\{ -t \, \widetilde{\Lambda}_i(n) \right\} \right].
		\end{split}
	\end{align}
	By \eqref{tilde} and the law of iterated expectation it follows that the expectation in the upper bound can be written as follows:  
	\begin{align} \label{markov}
		\Exp_A \left[    \exp \left\{ -t \, \widetilde{\Lambda}_i(n-1) \right\}
		\cdot  \Exp_A \left[  \exp\left\{-t \, (g_{i}(X_{i}(n)) -I_{i}) \, R_{i}(n) \right\} \, | \, \cF_{n-1}  \right]  \right].
	\end{align}
	Since $R(n)$ is  an $\cF_{n-1}$-measurable Bernoulli random variable and $X_i(n)$ is independent of $\cF_{n-1}$ and has the same distribution as $X_i(1)$, we also have 
	\begin{align*}
		&\Exp_{A}\left[\exp\left\{-t \, (g_{i}(X_{i}(n)) -I_{i}) \, R_{i}(n) \right\} \, | \, \cF_{n-1}  \right]  \\
		&=  \Exp_{A}\left[\exp\left\{-t \, (g_{i}(X_i(1)) -I_{i}) \right\}   \right]^{R_{i}(n)} \\
		&= \Exp_{A}\left[\exp\left\{t\,  (-g_{i}(X_i(1)) +I_{i} -\delta  + \delta) \right\}   \right]^{R_{i}(n)} \\
		& \leq  \exp\left\{ R_i(n) \, (\psi_{i}(t) + t\,  \delta  )\right\},
	\end{align*}
	where   $\psi_{i}$ is the cumulant generating function of $-g_{i}(X_{i}(1)) + I_{i} - \delta  $, i.e.,
	\begin{equation*}  
		\psi_{i}(t):=\log\Exp_{A}\left[\exp\left\{t (-g_{i}(X_i(1)) + I_{i} -\delta) \right\}   \right],  \;  t > 0.
	\end{equation*} 
	Since $\psi_i$ is convex on $(0,1)$, $\psi_{i}(1)=I_{i}-\delta <\infty$,  and  $$\psi_{i}'(0+)=\Exp_A[ -g_{i}(X_{i}(1)) + I_{i} - \delta ]=-\delta<0,$$ there is  an $s>0$ such that $\psi_{i}(s)<0$, and as a result 
	\begin{align*} 
		\exp\left\{ R_i(n) \, (\psi_{i}(s) + s\, \delta )\right\} &\leq  \exp\left\{s \, \delta  \right\}.
	\end{align*} 
	Therefore, setting $t=s$ in 
	\eqref{markov0}-\eqref{markov} we obtain  
	\begin{align*}
		\Pro_{A}  \left(  \widetilde{\Lambda}_i(n)<-n \, \epsilon/2 \right)
		&\leq  \exp\{-n \, (\epsilon/2)   s +\delta s  \}   \; \Exp_A\left[
		\exp \left\{ -s \, \widetilde{\Lambda}_i(n-1) \right\}  \right].
	\end{align*}
	Repeating the same argument $n-1$ times  we conclude that  there exists an $s>0$ such that 
	\begin{equation*}
		\Pro_{A}\left( \widetilde{\Lambda}_i(n)<-n \epsilon/2 \right)
		\leq \exp\left\{ -n (\epsilon/2-\delta) s   \right\}.
	\end{equation*}
	Since  $\delta \in (0, \epsilon/2)$, this  completes the proof. \\
\end{IEEEproof}

\begin{lemma}\label{new}
	Let $i\in A$,  $j \notin A$, $\zeta>0$, $\rho_i, \rho_j \in [0,1]$, and let $R$ be an arbitrary sampling rule.
	\begin{itemize}
		\item[(i)] If $\rho_i>0$ and $\Pro_{A}(\pi_i^R(n) <\rho_i )$ is summable, then so is $$\Pro_{A}\left(\frac{1}{n} \Lambda^R_i(n)< 
		\rho_i I_i- \zeta \right).$$
		\item[(ii)]  If  $\rho_j>0$ and $\Pro_{A} (\pi_j^R(n) <\rho_j )$ is summable, then so is $$\Pro_{A}\left(\frac{1}{n}  \Lambda^R_j(n)>-\rho_j J_j + \zeta \right).$$ 
		\item[(iii)]   If $\rho_i+\rho_j>0$ and both  $\Pro_{A}(\pi_j^R(n) <\rho_i )$ and 
		$\Pro_{A}(\pi_j^R(n) <\rho_j )$ are  summable, then so is 
		$$\Pro_{A}\left( \frac{1}{n} \Lambda^R_{ij}(n)
		< \rho_i I_i+ \rho_j J_j -\zeta \right).
		$$
	\end{itemize}
\end{lemma}

\begin{IEEEproof} For (i), by the law of total probability  we have
	\begin{align*}
		&\Pro_{A} \left( \frac{1}{n} \Lambda^R_i(n)< \rho_i I_i-\zeta \right) \\
		&\leq  \Pro_{A}\left(\pi_i^R(n) <\rho_i \right)  + \Pro_{A}\left(\frac{1}{n} \Lambda^R_i(n)< 
		\rho_i I_i-\zeta, \,  \pi^R_i(n) \geq \rho_i \right).
	\end{align*}
	The first term in the upper bound is summable by assumption, whereas the second one is summable, as exponentially decaying,   by  Lemma \ref{lem:A1}. 
	
	The proof of (ii) is similar and is omitted.  Since  (iii) follows directly by (i) and (ii) when $\rho_i$ and $\rho_j$ are both positive,  it suffices to consider the case where only one them  is positive. Without loss of generality, we assume that $\rho_i>0=\rho_j$. Then,  by the law of total probability again we have
	\begin{align*}
		&\Pro_{A}\left( \frac{1}{n}  \Lambda^R_{ij}(n)  < \rho_i I_i+ \rho_j J_j -\zeta  \right)
		\\
		& \leq  \Pro_{A}\left(\pi_i^R(n) <\rho_i \right) + \Pro_{A}\left(\frac{1}{n} \Lambda^R_{ij}(n)  < 
		\rho_i I_i-\zeta, \,  \pi^R_i(n) \geq \rho_i \right).
	\end{align*}
	The first term in the upper bound is summable by assumption, whereas the second one is  summable, as exponentially decaying, by  Lemma \ref{lem:A1}.  \\
\end{IEEEproof}

\section{} \label{app:consistency_criterion}

In this Appendix we prove Theorems \ref{propo_1} and  \ref{theo:prob_cons},  which provide sufficient conditions for the exponential consistency of a sampling rule. In order to lighten the notation, throughout this Appendix   we suppress  dependence on $R$ and  we write 
$\pi_i(n)$,  $\Lambda_i(n), \Lambda_{ij}(n),\Delta_n$, $c_i(A)$, $\sigma_A$ instead of
$\pi_i^R(n)$, $\Lambda_i^R(n)$, $\Lambda_{ij}^R(n)$, $\Delta^R_n$, $c^R_i(A)$, $\sigma^R_A$.


\begin{IEEEproof}[Proof of Theorem \ref{propo_1}]
	We prove the result first when $\ell=u$. By the definition of the decision rule in \eqref{gap_decision rule} it follows that when  $\sigma_A> n$, there are  $m \geq  n$, $i \in A$, $j \notin A$ such that  $\Lambda_{ij}(m)\leq 0$,
	and as a result  by the union bound we have 
	\begin{equation*}
		\Pro_{A}(\sigma_{A} > n) \leq   \sum\limits_{i \in A, \, j \notin A} \; \sum_{m= n }^\infty  \Pro_{A} \left(\Lambda _{ij}(m)\leq  0 \right). 
	\end{equation*}
	Therefore, to show that $\Pro_{A}(\sigma_{A} > n)$ is exponentially decaying, it suffices to show that this is the case for $\Pro_{A} \left(\Lambda _{ij}(n)\leq  0 \right)$ for every  $i \in A$ and  $j \notin A$.  To show this, we fix such $i$ and $j$ and we note that, by assumption, either $\Pro_{A}(\pi_{i}(n)< \rho)$  or $\Pro_{A}(\pi_{j}(n)< \rho)$ is exponentially decaying  for  $\rho>0$ small enough. Without loss of generality, suppose that this is the case for  the former. By an application of the law of total probability  we then  obtain 
	\begin{align*}
		\Pro_{A} \left( \Lambda _{ij}(n) \leq 0 \right)
		&\leq \Pro_{A}(\Lambda _{ij}(n) \leq 0,  \pi_i(n) \geq \rho) \\
		&+ \Pro_A(\pi_i(n) < \rho).
	\end{align*}	
	As mentioned earlier, the second term in the upper bound is, by assumption,  exponentially decaying for  $\rho>0$ small enough. By  Lemma \ref{lem:A1} it follows that this is also the  case for  the first one.
	
	Suppose next that  $\ell<|A|<u$. By the definition of the decision rule in \eqref{gi_decision rule} it follows that when  $\sigma_A>  n$,  there is an $m \geq  n$ such that either  $\Lambda_{i}(m)<0$ for some   $i \in A$, or  $\Lambda_{j}(m)>0$ for some $j \notin A$. As a result, by the union bound we have 
	\begin{align*}
		\Pro_{A}({\sigma}_{A} >  n) 
		&\leq \sum\limits_{i \in A}   \; \sum_{m= n }^\infty \Pro_{A}(\Lambda_{i}(m)<0) \\
		&+ \sum\limits_{j \notin A}  \; \sum_{m= n  }^\infty \Pro_{A}( \Lambda _{j}(m) \geq 0).
	\end{align*}
	Therefore, to prove that $\Pro_{A}(\sigma_{A} >  n)$ is exponentially decaying, it suffices to show that this is true for $\Pro_{A} \left(\Lambda _{i}(n) < 0 \right)$ and   $\Pro_{A} \left(\Lambda _{j}(n)\geq  0 \right)$  for every  $i \in A$ and   $j \notin A$.   This can be shown similarly to the case $\ell=u$, in  view of the fact that  the sequences 
	$\Pro_{A}(\pi_{i}(n)< \rho)$ and $\Pro_{A}(\pi_{j}(n)< \rho)$  are,  by assumption, both exponentially decaying  for $\rho>0$ small enough and  for every  $i \in A$ and   $j \notin A$. 
	
	The two remaining cases are  $\ell<|A|=u$ and $\ell = |A|<u$. We consider only the former, as the proof for the latter is similar, and assume  that  $\ell<|A|=u$. By the definition of the decision rule in \eqref{gap_decision rule} it follows that when  $\sigma_A>  n$, there is an $m\geq  n$ such that  either $\Lambda_{ij} (m) \leq 0$ for some $i \in A$ and $j \notin A$, or $ \Lambda_{i}(m)<0$ for some $i \in A$. As a result, by the union bound we have 
	\begin{align*}
		\Pro_{A}(\sigma_{A} > \zeta n) &\leq  \sum\limits_{i \in A}   \sum_{m= n }^\infty \Pro_{A}(\Lambda_{i}(m)<0) \\
		&+	  \sum\limits_{i \in A,j \notin A}   \sum_{m= n}^\infty \Pro_{A}(\Lambda _{ij}(m) \leq 0).
	\end{align*}
	Since  $|A|>\ell$,  the sequence $\Pro_{A}(\pi_{i}(n)< \rho)$ is, by assumption,  exponentially decaying for every $i \in A$ and $\rho>0$ small enough.  Therefore, similarly  to the previous cases we can show that each term in the upper bound is exponentially decaying.\\
\end{IEEEproof}


In the remainder of this Appendix we prove Theorem \ref{theo:prob_cons}, whose proof relies on two preliminary lemmas. To state those,   for every $D \in \mathcal{P}_{\ell, u}$ and  $n \in \bN$  we denote by  $\tau^D_n$ the first time instance $m$  at which  $D$ has been estimated as the subset of anomalous sources  for at least  $2 \, \zeta \,  m$ times since   $ \lceil n/2  \rceil $, i.e.  
\begin{equation}\label{tau}
	\tau^D_{n}:=\inf\left\{ m \geq \lceil n/2 \rceil : \sum_{u=\lceil n/2\rceil}^{m} \mathbf{1}\{\Delta_{u-1}=D\}  \geq 2 \zeta    m \right\},
\end{equation}
where  $\zeta>0$ is an  arbitrary constant, which in the proof of Theorem \ref{theo:prob_cons} will be selected to be small enough.  \\

\begin{lemma}\label{lem:prob_consi_lem0}
	Let  $D \in \mathcal{P}_{\ell, u}$ and $i \in [M]$. If   $c_{i}(D)>0$, 	then
	$\Pro(\tau^D_{n}  \leq  n,\;  \pi_{i}(\tau^D_n)<\rho)$ and 
	$\Pro(\tau^D_{n} \leq  n,\;  \pi_{i}(n)<\rho)$ are both  exponentially decaying for all $\rho>0$ small enough.  \\
\end{lemma}

\begin{IEEEproof} 
	We prove the two claims together by showing that $\Pro(\tau^D_{n} \leq  n,\;  \pi_{i}(\sigma_n)<\rho)$ is exponentially decaying  for all $\rho>0$ small enough, where $\sigma_n$ stands for either $n$ or $\tau_n^D$. For any given $n \in \bN$ we set 
	\begin{align*}
		\widetilde{\pi}_{i}(n) &:=   \frac{1}{n} \sum_{m=1}^{n} (R_i(m)- c_{i}(\Delta_{m-1}) ) \\
		&= \pi_i(n) -\frac{1}{n} \sum_{m=1}^{n} c_{i}(\Delta_{m-1}), \quad n \in \bN.
	\end{align*}  
	Then,  on the event  $\{  \tau^D_n \leq n\}$ we have 
	$n/2\leq \tau^D_n \leq \sigma_n \leq n$ and consequently
	\begin{align*}
		\pi_i(\sigma_n) -\widetilde{\pi}_i(\sigma_n)   &= \frac{1}{\sigma_n}  \sum_{m=1}^{\sigma_n} c_{i}(\Delta_{m-1}) \\
		&\geq   \frac{c_i(D)}{n} \sum_{m=\lceil n/2 \rceil}^{\tau^D_n}   \, \mathbf{1}\{\Delta_{m-1}=D\}\\
		&\geq \frac{c_i(D)}{n} \, 2 \, \zeta \,  \tau_n^D \geq \; c_i(D) \, \zeta,
	\end{align*}
	where in the last inequality we have used  the definition of $\tau_n^D$. Consequently,  for every $\rho>0$ and $n \in \bN$   we obtain
	$$\{\tau^D_{n} \leq  n,\;  \pi_{i}(\sigma_n) \leq \rho \} \subseteq  \{ \tau^D_{n} \leq  n, \, \widetilde{\pi}_i(\sigma_n) \leq \rho-c_{i}(D) \, \zeta \}.
	$$
	Since $c_i(D)>0$,  there is an $\epsilon>0$ such that 
	for all  $\rho \in (0, c_{i}(D) \zeta -\epsilon)$ we  have 
	$$\{\tau^D_{n} \leq  n,\;  \pi_{i}(\sigma_n) \leq \rho \}  \subseteq  
	\left\{ \max_{\lceil n/2 \rceil \leq m \leq n} |\widetilde{\pi}_i(m)| 
	\geq  \epsilon \right\}.
	$$
	Therefore, it suffices to show that,  for all $\epsilon>0$, the sequence
	\begin{align}\label{azuma}
		\Pro\left(\max_{\lceil n/2 \rceil \leq m \leq n} |\widetilde{\pi}_i(m)| 
		\geq  \epsilon \right) 
	\end{align}
	is exponentially decaying. This follows by the fact that  $(R_{i}(n) -c_{i}(\Delta_{n-1}))$ is a uniformly bounded  martingale difference, in view of   \eqref{c}, and an application of the  maximal Azuma-Hoeffding submartingale inequality \cite[Remark 2.3]{rag_sas13}. \\
\end{IEEEproof}

\begin{lemma}\label{lem:prob_consi_lem1}
	Let  $A, D \in \mathcal{P}_{\ell, u}$  be such that  $|D| \in \{l, u\}$ and $D \neq A$. If  $R$  is a sampling rule that satisfies the conditions of Theorem \ref{theo:prob_cons}, then  the sequence $\Pro_{A}(\tau^D_{n} \leq n)$ is exponentially decaying.\\
\end{lemma}

\begin{IEEEproof}  We consider first the case that $\ell<u$, where we only prove the result when $|D|=u$, as the proof when $|D|=\ell$ is similar. Since  $A \in \mathcal{P}_{\ell, u}$,  $|D|=u$ and $D \neq A$, there exists a $j \in D \setminus A$. Since also $|D|>\ell$, the assumptions of Theorem \ref{theo:prob_cons}  imply that   $c_j(D)>0$. Thus, by  Lemma \ref{lem:prob_consi_lem0}  it follows that
	\begin{align} \label{pr0}
		\Pro_{A}\left( \tau^D_{n} \leq n, \; \pi_{j}(\tau^D_{n})<\rho \right) 
	\end{align}
	is an exponentially decaying  sequence for $\rho>0$ small enough,  and  it suffices to show that this is also the case for 
	\begin{align} \label{pr1}
		\Pro_{A}\left( \tau^D_{n} \leq n,  \; \; \pi_{j}(\tau^D_{n}) \geq \rho  \right) .
	\end{align}
	Indeed,  by the  definition  of $\tau^D_{n}$ in  \eqref{tau}  it follows that  on the event $\{\tau^D_{n} <\infty \}$ we have $\Delta(\tau_{n}^D-1)=D$. Since  $|D|=u$ and $j \in D$, by the definition of the decision rule in  \eqref{gi_decision rule}  it follows that $\Lambda_{j}(\tau^D_{n}-1)>0$.  Therefore, 
	\begin{align*} 
		\begin{split}
			&\Pro_{A}\left( \tau^D_{n} \leq n,  \; \; \pi_{j}(\tau^D_{n})  \geq \rho \right) \\
			&= \Pro_{A}( \tau^D_n\leq n, \,  \Lambda_{j}(\tau_{n}^D-1)>0,\,  \pi_{j}( \tau_n^D) \geq \rho) \\
			&\leq   \sum_{m=\lceil n/2  \rceil }^n \Pro_{A}\left( \Lambda_{j}(m-1) >0, \, \pi_{j}(m) \geq  \rho \right) ,
		\end{split}
	\end{align*}
	where the  inequality holds because  $\tau_n^D$ takes values in $[n/2, n]$ on the event $\{\tau_n^D \leq n\}$.  Therefore, it remains to show that the sequence
	\begin{align} \label{exist_11} 
		\Pro_{A}\left( \Lambda_{j}(n-1) >0, \, \pi_{j}(n) \geq  \rho \right)
	\end{align}
	is exponentially decaying for $\rho>0$ small enough. Indeed, for large enough $n$,  $\pi_{j}(n) \geq \rho$ implies that 
	$$\pi_{j}(n-1) \geq \frac{n\pi_{j}(n)-1}{n-1} \geq \frac{\rho n -1}{n-1} \geq \rho - \frac{1}{n-1}.$$
	For any given $ \rho>0$, there exists a $\rho'>0$ so that for all $n$ large enough we have 
	\begin{align} \label{rho_prime}
		\rho - \frac{1}{n-1} > \rho',
	\end{align}
	and so that the  probability in \eqref{exist_11} is bounded by 
	$$\Pro_{A}\left( \Lambda_{j}(n-1) >0, \, \pi_{j}(n-1) \geq  \rho' \right).$$
	For $\rho>0$ small enough,  $\rho'>0$ is  small enough,  and by Lemma \ref{lem:A1} it follows that the latter probability, and consequently  \eqref{exist_11}, is exponentially decaying in $n$.
	

	
	It remains to prove the lemma when $\ell=u$ and $A \neq D$.  In this case there  are $i \in A \setminus D$ and $j \in D \setminus A$, and by the assumptions of Theorem \ref{theo:prob_cons} it follows that either 
	$c_i(D)>0$ or $c_j(D)>0$.   Without loss of generality, we assume that the latter holds. Then, by Lemma \ref{lem:prob_consi_lem0}  it  follows that  \eqref{pr0}  is exponentially decaying for all $\rho>0$ small enough, and it  suffices to show that  this is also the case for  \eqref{pr1}.  Indeed, by the  definition  of $\tau^D_{n}$ in  \eqref{tau}  it follows that  on the event $\{\tau^D_{n}<\infty\}$ we have $\Delta(\tau_{n}^D-1)=D$.  Consequently, by the definition of the decision rule in  \eqref{gap_decision rule}
	it follows that there is an $i \in A \setminus D$ such that 
	$\Lambda_{ij}(\tau^D_{n}-1)<0$.  Therefore,  by the union bound we have 
	\begin{align*}
		&\Pro_{A}\left( \tau^D_{n} \leq n,  \; \; \pi_{j}(\tau^D_{n})  \geq \rho \right) \nonumber \\
		& = \sum_{i \in A \setminus D} \Pro_{A} \left( \tau^D_n\leq n, \,  \Lambda_{ij}(\tau^D_{n}-1) <0,\,  \pi_{j}( \tau_n^D) \geq \rho \right) \nonumber \\
		&\leq \sum_{i \in A \setminus D}  \sum_{m=\lceil n/2  \rceil}^n  \Pro_{A}\left(\Lambda_{ij}(m-1)<0, \, \pi_{j}(m) \geq  \rho \right),
		\label{exist_2}
	\end{align*}
	where as before the  inequality holds because  $\tau_n^D$ takes values in $[n/2, n]$ on the event $\{\tau_n^D \leq n\}$. Therefore, it remains to show that the sequence
	\begin{align*} 
		\Pro_{A}\left( \Lambda_{ij}(n-1) <0, \, \pi_{j}(n) \geq  \rho \right)
	\end{align*}
	is exponentially decaying for $\rho>0$ small enough. As before, this follows by an application of Lemma  \ref{lem:A1}.
\end{IEEEproof}

\begin{IEEEproof}[Proof of Theorem \ref{theo:prob_cons}]   
	Fix $A \in \mathcal{P}_{\ell, u}$.   By Theorem  \ref{propo_1}  it suffices to show that, for all $\rho>0$ small enough,  $\Pro_A(\pi_i(n) < \rho)$  is an exponentially decaying sequence 
	\begin{itemize} 
		\item for every $i \in A$, if  $|A|>\ell$, and  for every $i \notin A$, if  $|A|<u$, when $\ell<u$,
		\item  either for every $i \in A$ or for every $i \notin A$, when  
		$\ell=u$. 
	\end{itemize} 
	In order to do so, we select  the positive constant $\zeta$ in \eqref{tau} to be  smaller than  $1 / |\mathcal{P}_{\ell, u}|$. Then, for every  $n \in \bN$   there is at least one $D \in \mathcal{P}_{\ell, u}$ for which $\{\tau^D_n \leq n \} \neq \emptyset$. As a result, for every   $i \in [M]$ and $\rho>0$, by the union bound   we have 
	\begin{equation*}
		\Pro_{A}(\pi_{i}(n)<\rho) \leq \sum_{D\in \mathcal{P}_{\ell, u}} \Pro_{A}\left( \pi_{i}(n)<\rho, \; \tau_n^D\leq n \right).
	\end{equation*}
	Suppose first that $\ell < u$.  Then,  it suffices to show that, for every $D \in \mathcal{P}_{\ell, u}$  and all $\rho>0$ small enough,  
	\begin{equation}\label{boo_f}
		\Pro_A\left(\pi_i(n) < \rho, \tau_n^D\leq n \right)
	\end{equation} 
	is exponentially decaying  for every $i \in A$   when $|A|>\ell$ and for  every  $i \notin A$ when  $|A|<u$.   We  only consider  the former case, as the proof for the latter is similar. Thus, suppose that $|A|>\ell$ and let  $i \in A$.
	\begin{itemize}
		\item    If $c_{i}(D)>0$,   by Lemma \ref{lem:prob_consi_lem0}  it follows that \eqref{boo_f}  is an exponentially decaying sequence.
		\item  If  $c_{i}(D)=0$, the assumption of the theorem implies that   either    $|D|=u$ and $ i \notin D$, or $|D|=\ell$ and $i \in D$.  
		In either case,  $ A \neq  D$  and 
		by Lemma \ref{lem:prob_consi_lem1} it follows that $\Pro_{A}(\tau_n^D\leq n )$, and consequently \eqref{boo_f},  is an exponentially decaying sequence.
	\end{itemize}

	Suppose now that    $\ell=u$.  Then,  it suffices to  show that \eqref{boo_f}   is exponentially decaying  for every $D \in \mathcal{P}_{\ell, u}$ and  all $\rho>0$ small enough,   either for every $i \in A$ or for every $i \notin A$. 
	
	\begin{itemize}
		\item When $D \neq A$, by Lemma  \ref{lem:prob_consi_lem1} it follows that $\Pro_{A}(\tau_n^D\leq n )$, and consequently \eqref{boo_f}, is exponentially decaying.
		\item  When $D=A$, then by  assumption $c_{i}(A)>0$ holds  either  for every $i \in A$ or  for every  $i \notin A$. By   Lemma  \ref{lem:prob_consi_lem0} it then follows  that \eqref{boo_f}  is exponentially decaying either for every  $i \in A$ or for every  $i \notin A$, and this completes the proof. 
	\end{itemize}
\end{IEEEproof}

\section{} \label{app:lower_bound} 
In this  Appendix we fix $A \in \cP_{\ell, u}$ and  prove the universal asymptotic lower bound of Theorem \ref{theo:AO}.  The proof relies on two lemmas, for the statement of which we need to  introduce the following function:
\begin{align*}
	\phi(\alpha,\beta) := \alpha \log\left(\frac{\alpha}{1-\beta} \right) + (1-\alpha) \log\left(\frac{1-\alpha}{\beta} \right), 
\end{align*}
where $\alpha, \beta \in (0,1)$,  i.e., the Kullback-Leibler divergence between a Bernoulli distribution with parameter $\alpha$  and one with parameter $1-\beta$. Moreover, we set $\phi(\alpha) \equiv  \phi(\alpha,\alpha)$.  

The first lemma states a non-asymptotic, information-theoretic inequality that   generalizes  the one used in Wald's universal lower bound in the problem of testing two simple hypotheses \cite[p. 156]{wald1945sequential}.\\

\begin{lemma} 
	Let  $\alpha+\beta <1$ and let $(R, T, \Delta)$ be a policy  that satisfies the error constraint \eqref{err_const} and  $\Pro_A(T<\infty)=1$.  Then, for any  $C \in \cP_{\ell, u}$ such that $C \neq A$ we have 
	\begin{equation} \label{info_ineq}
		\Exp_{A} \left[ \Lambda^R_{A,{C}}(T) \right] \geq 
		\begin{cases}
			\phi(\alpha,\beta) \; & \text{if } \;  {C}{\setminus}A \neq \emptyset,\mbox{ } A{\setminus}{C} = \emptyset , \\
			\phi(\beta,\alpha) \; & \text{if } \; {C}{\setminus}A = \emptyset,\mbox{ } A{\setminus}{C} \neq \emptyset, \\
			\phi(\alpha \wedge \beta) \quad &\text{if } \;  {C}{\setminus}A \neq \emptyset,\mbox{ } A{\setminus}{C} \neq \emptyset. 
		\end{cases}
	\end{equation}	
\end{lemma}

\begin{IEEEproof}
	The proof is identical to that in the full sampling case
	in \cite[Theorem 5.1]{Song_and_Fellouris_2016}, and  
	can be obtained by an application of the data processing inequality for Kullback-Leibler divergences  (see, e.g., \cite[Lemma 3.2.1]{tartakovsky_book_2014}). Indeed,   the left-hand side is the  Kullback-Leibler  divergence between $\Pro_A$ and $\Pro_C$ given the available information up  to time $T$, when the sampling rule $R$ is utilized, whereas the right hand side is obtained by  considering the Kullback-Leibler divergence between $\Pro_A$ and  $\Pro_C$ based on  a single event of  $\mathcal{F}^R_T$.
\end{IEEEproof}

We next make use of the previous  inequality to establish lower bounds on the expected number of samples taken from each source until stopping. \\

\begin{lemma} \label{lem:appenB3}
	Let  $\alpha+\beta <1$ and let $(R, T, \Delta)$ be a policy  that satisfies the error constraint \eqref{err_const} and  $\Exp_A[T]<\infty$. 
	\begin{enumerate}
		\item[(i)]  If  $|A|<u $, then 
		\begin{equation}\label{take_exp_notA}
			\min_{j \notin A}  \left( J_j \, \Exp_{A}\left[ N_j^R(T) \right] \right)  \geq \phi(\alpha, \beta).
		\end{equation} 
		\item[(ii)] If $|A|>\ell $, then 
		\begin{equation}\label{take_exp_A}
			\min_{i \in A}  \left( I_i \, \Exp_{A}\left[ N_i^R(T)  \right] \right)  \geq \phi( \beta, \alpha) .
		\end{equation} 
		\item[(iii)] 	If either $|A|=\ell>0$ or $|A|= u<M$, then 
		\begin{align}  \label{take_exp_both}
			&\min_{i \in A} \left( I_i \, \Exp_A\left[ N_i^R(T)  \right]\right) 
			+\min_{j \notin A}	 \left( J_{j} \, \Exp_A\left[ N_j^R(T) \right] \right)   \geq \phi(\alpha \wedge \beta) .
		\end{align}
	\end{enumerate}
\end{lemma}

\begin{IEEEproof} 
	We recall the  sequence  
	$\widetilde{\Lambda}^R_i,$  defined in \eqref{tilde}, and note  that   it is a  zero-mean,  $\{ \cF^R_{n} \}$-martingale under  $\Pro_A$. Moreover,  by the finiteness of the Kullback-Leibler divergences in \eqref{KL}  we have:
	\begin{equation*}  
		\sup_{n \in \bN}	\Exp_A\left[ | \widetilde{\Lambda}^R_i(n)-\widetilde{\Lambda}^R_i(n-1)| \;  | \;  \cF^R_{n-1} \right] <\infty.
	\end{equation*}
	Since also   $T$ is an $\{\mathcal{F}^R_n\}$-stopping time such that $\Exp_A[T]<\infty$,  by  the  Optional Sampling Theorem \cite[pg. 251]{Chow_Teicher2012} we obtain:
	\begin{equation} \label{OST}
		\Exp_A\left[ \widetilde{\Lambda}^R_i(T) \right]=0 \quad \text{for every}  \quad i \in [M].
	\end{equation}
	
	(i) If  $|A|<u $, there is a $j \notin A$ such that   $C=A \cup \{j\} \in \cP_{\ell, u}$.  By representation \eqref{LLR_global_repre} and 
	decomposition  \eqref{decompose} it follows that   the log-likelihood ratio process in \eqref{LLR_global} takes the form 
	\begin{equation*}
		\Lambda_{A,C}^R(T)=-\Lambda_{j}^R(T)=-\widetilde{\Lambda}_{j}^R(T)+J_{j} \, N_j^R(T) .
	\end{equation*}
	By  \eqref{info_ineq} and   \eqref{OST}  we then obtain 
	$$J_{j} \, \Exp_{A}\left[ N_j^R(T) \right] \geq \phi(\alpha,\beta).$$ 
	Since this inequality holds for every $j \notin A$, this  completes the proof. 
	
	(ii) The proof is similar to (i) and is omitted. 
	
	(iii) If  $|A|=\ell>0$ or $|A|= u<M$, then there are  $i \in  A$ and  $j \notin A$ such that   $C=A \cup \{j\} \setminus \{i\} \in \cP_{\ell, u}$. By  representation \eqref{LLR_global_repre} and  decomposition \eqref{decompose} we have 
	\begin{align*}
		\Lambda^R_{A,C}(T) &= \Lambda_i^R(T)- \Lambda_{j}^R(T) \\
		&=\widetilde{\Lambda}_i^R(T)- \widetilde{\Lambda}^R_j(T)+J_{j}  \, N_j^R(T) + I_{i} \, N_i^R(T) .
	\end{align*}	
	By  \eqref{info_ineq} and   \eqref{OST}  we then obtain
	\begin{equation*}
		I_i \,  \Exp_A\left[ N_i^R(T) \right]+	J_{j} \, \Exp_A\left[ N_j^R(T)  \right]  \geq  \phi(\alpha \wedge \beta).
	\end{equation*}
	Since this inequality holds  for every $i \in A$ and $j \notin A$, this proves \eqref{take_exp_both}.\\
\end{IEEEproof}

For the proof of Theorem \ref{theo:lowerbound} we  introduce the following notation:
\begin{align}\label{domain}
	\begin{split}
		\cD_{K} &:=\left\{ (c_1, \ldots, c_M) \in [0,1]^M:  \sum_{i=1}^M  c_i\leq K\right\} \\
		\cD'_K &:=\{(p,q) \in [0,1]^2: \; p \hat{K}_A+q \check{K}_A \leq K\}.
	\end{split}
\end{align} 
Moreover, we  observe that as $\alpha, \beta \to 0$ we have  
\begin{align}
	\phi(\alpha, \beta) & \sim  |\log \beta|,  \label{phi}\\
	r(\alpha, \beta) \equiv \frac{\phi(\beta, \alpha)} {\phi(\alpha, \beta)} &\sim 
	\frac{|\log \alpha|}{|\log \beta|}.
	\label{phi_ratio}
\end{align} 


\begin{IEEEproof}[Proof of Theorem \ref{theo:lowerbound}] 
	(i) Let  $\alpha, \beta \in (0,1)$ such that $\alpha+\beta<1$  and    $(R,T,\Delta) \in \cC(\alpha, \beta, \ell, u,K)$ such that $\Exp_A[T]<\infty$.  By  Lemma  \ref{lem:appenB3}(iii) it then follows that 
	\begin{align*}
		&\Exp_A[T] \;  W_A(T) \geq   \phi(\alpha \wedge \beta), \quad 
		\text{where} \\
		W_A(T) &:=
		\min_{i \in A}  \left\{ I_{i} \,  \frac{\Exp_A[N^R_i[T]}{\Exp_A[T]} \right\}+\min_{j \notin A} \left\{ J_j \,  \frac{\Exp_A[N^R_j[T]}{\Exp_A[T]}  \right\},
	\end{align*}
	and by constraint \eqref{samp_const} we conclude that
	\begin{align*}
		\Exp_A[T] \;  V_A &\geq   \phi(\alpha \wedge \beta), \\
		\text{where} \qquad 	V_A &:= \max_{(c_1, \ldots, c_M) \in \mathcal{D}_K} \left\{ 	\min\limits_{i \in A} (c_{i}I_{i})+\min\limits_{j{\notin}A}(c_{j} J_{j}) \right\} .
	\end{align*}
	Since  the lower bound  is independent of the policy $(R,T,D)$, we have 
	\begin{align*}
		\cJ_A(\alpha, \beta,\ell, u, K)\;  V_A &\geq   \phi(\alpha \wedge \beta).
	\end{align*}
	Comparing with  \eqref{LB_gap} 
	and recalling \eqref{phi}, we can see that it suffices to show that 	$V_A=x_A  I_A^* +y_A J_A^*$  with $x_A$ and $y_A$ as in \eqref{xy_gap1}-\eqref{xy_gap}. 
	Indeed, the maximum in $V_A$ is  achieved by $c_i$'s  of the form 
	\begin{align} \label{optimizers}
		\begin{split}
			c_i I_i &= p I_A^*, \quad i \in A,\\
			c_j J_j &= q J_A^*, \quad j \notin A,
		\end{split}
	\end{align}
	for  $p,q \in [0,1]$  such that the constraint in $\cD_K$ is satisfied, i.e., 
	\begin{equation} \label{constraint}
		K\geq \sum_{i=1}^M c_i = p \sum_{i \in A}\frac{ I^{*}_{A}} 
		{I_{i}}  +q \sum_{j \notin A}  \frac{ J^{*}_{A}} 
		{J_{j}} = p \, \hat{K}_A+ q\, \check{K}_A,
	\end{equation}
	and as a result, 	
	$$V_A =\max_{(p,q) \in \cD'_K}  \{p I_A^*+ q J_A^*\}.$$
	This maximum is achieved by  $p,q \in [0,1]$ such that 
	$p \hat{K}_A+ q \check{K}_A = K \wedge (\hat{K}_A+\check{K}_A)$, in particular by  $p$ and $q$ equal to $x_A$ and $y_A$  as in  \eqref{xy_gap1}-\eqref{xy_gap}, which completes the proof. \\
	

	(ii)  Suppose first that $ \ell<|A|<u $.
	As before, let  $\alpha, \beta \in (0,1)$ such that $\alpha+\beta<1$  and    $(R,T,\Delta) \in \cC(\alpha, \beta, \ell, u,K)$ such that $\Exp_A[T]<\infty$.   Then, by  
	Lemma \ref{lem:appenB3}(i) 
	and Lemma \ref{lem:appenB3}(ii) we obtain:
	\begin{align*}
		&\Exp_A[T] \;  W_A(T) \geq   \phi( \beta, \alpha ), \;
		\text{where} \\
		W_A(T) &:=  \min \left\{
		\min_{i \in A}  \left\{ I_{i} \,  \frac{\Exp_A[N^R_i[T]}{\Exp_A[T]} \right\} , \;  r(\alpha, \beta) \; \min_{j \notin A} \left\{ J_j \,  \frac{\Exp_A[N^R_j[T]}{\Exp_A[T]}  \right\}  \right\},
	\end{align*}
	and  by constraint \eqref{samp_const} we conclude that
	\begin{align}	\label{V_int}
		\begin{split}
			& \Exp_A[T] \;  V_{A}(\alpha, \beta) \geq 
			\phi( \beta, \alpha )  \\
			\text{where} \qquad 	V_A (\alpha, \beta)  &:= \max_{(c_1,\ldots, c_M) \in \cD_K} \min\left   \{ \min_{i \in A} \left(c_i I_i \right), \, r(\alpha, \beta)   \;  \min_{j \notin A}  \left(c_j J_j\right)\right\}.
		\end{split}
	\end{align}
	Since the lower bound does not depend on  the policy $(R,T,\Delta)$  we have 
	\begin{align*}
		\cJ_A(\alpha, \beta,\ell, u, K) \;  V_{A}(\alpha, \beta) &\geq 
		\phi( \beta, \alpha ).
	\end{align*}
	Comparing with	\eqref{LB_int}  and recalling  \eqref{phi} we can  see that it suffices to show that 
	$$V_A(\alpha, \beta)  \rightarrow  x_A \,  I_A^* =  r
	\, y_A \,  J_A^*$$
	as 
	$\alpha, \beta \to 0$ according to \eqref{r},   with  $x_A$ and $y_A$  as in  \eqref{LB_int}. 
	The equality follows directly from the values of  $x_A$ and $y_A$ in \eqref{LB_gap_int_1}. Moreover, 
	the  maximum in $V_A(\alpha, \beta)$   is achieved by $c_1, \ldots, c_M$ of the form \eqref{optimizers} that satisfy  \eqref{constraint}. Therefore:
	$$V_A(\alpha, \beta) =\max_{(p,q) \in \cD'_K}  \min\left\{ p I_A^*, \,   r(\alpha, \beta)  \, q J_A^* \right\},$$
	and this   maximum is achieved for $p$ and $q$ such that the two terms in the minimum are equal.  As a result, we obtain 
	$$V_A(\alpha, \beta) = p I_A^* = r(\alpha, \beta)  \, q J_A^*,$$
	where  $p$  and $q$  are equal to $x_A$ and  $y_A$  in  \eqref{LB_int},  with $r$  replaced by 
	$r(\alpha, \beta)$.   As  $\alpha$ and  $\beta$ go to  $0$  according to \eqref{r}, we have $r(\alpha, \beta) \to r$  (recall \eqref{phi_ratio}),  and consequently  $V_A(\alpha, \beta)  \rightarrow  x_A \,  I_A^*$   with $x_A$ as  in \eqref{LB_gap_int_1}.\\



	
	Finally, we consider the case   $|A|=\ell$ and omit the  proof when $|A|=u$, as it is similar.  When  either $\ell=0$ or $r \leq 1$, we have to show  \eqref{LB_00}.  Indeed, working as before, using Lemma \ref{lem:appenB3}(i), we obtain 
	\begin{align*} 
		\cJ_A(\alpha, \beta,\ell, u, K) \;  V_{A}  &\geq \phi(\alpha, \beta), \\
		\text{where} 	 \quad 	V_A &:= \max_{(c_1, \ldots, c_M) \in \cD_K}  \; \min_{j \notin A}\; \left\{ c_j J_j \right\} .
	\end{align*}
	Comparing with	\eqref{LB_00}, and recalling  \eqref{phi},  we can  see that it suffices to show that 
	$	V_A=  J_A^*\;   y_A,
	$ with $y_A$ as in \eqref{LB_00}. Indeed, the  maximum in   $V_A$    is achieved by  
	$c_1, \ldots, c_M$  of the form \eqref{optimizers} with $p=0$ and  
	$q \in [0,1]$ such that   \eqref{constraint} is satisfied, i.e.,
	$$	V_A = J_A^*\;   \max_{q \in [0,1] : \; q \check{K}_A \leq K} q,
	$$
	which shows that 
	$	V_A=  J_A^*\;   y_A,
	$ with $y_A$  as in \eqref{LB_00} and completes the proof in this case.

	It remains to establish the asymptotic lower bound when  $\ell>0$ and $r > 1$,  in which case  we have to show  \eqref{LB_gap_int_1}-\eqref{LB_gap_int_2}.  The asymptotic equivalence in \eqref{LB_gap_int_1} can be shown by direct evaluation,  therefore it suffices to show only the asymptotic lower bound  in this case. 
	
	Working as before, using  Lemma \ref{lem:appenB3}(i) and 
	Lemma \ref{lem:appenB3}(iii), for any $\alpha,  \beta \in (0,1)$ such that $\alpha +\beta <1$  we obtain
	\begin{align}  
		\label{V_A_bou1} 
		\begin{split}
			&\cJ_A(\alpha, \beta,\ell, u, K) \; \; V_A(\alpha, \beta) \geq  \phi(\beta, \alpha),  \quad
			\text{where} \\
			V_A(\alpha, \beta)  &:= \max_{(c_1, \ldots, c_M) \in \cD_K} \min\left\{r(\alpha, \beta) \; \min_{j \notin A}\left(c_j J_j\right),  \, \min_{i{\in}A}\left(c_i I_i\right)+\min_{j \notin A}\left(c_j J_j\right)  \right\} .
		\end{split}
	\end{align} 
	The latter maximum  is achieved by $c_1, \ldots, c_M$ of the form \eqref{optimizers} that satisfy \eqref{constraint}, thus,
	\begin{align*}
		V_A(\alpha, \beta)  &=	\max_{(p,q) \in \cD'_K}  \; \min\left\{r(\alpha, \beta) \; q J_A^*,\;  p I_A^*+q J_A^*\right\}.
	\end{align*}

	\begin{itemize}
		
		\item If   $ \theta_A \geq r(\alpha, \beta)-1$ or  $K\leq  \hat{K}_A+  (\theta_A/(r(\alpha, \beta)-1)) \check{K}_A$,  the maximum in $V_A(\alpha, \beta)$  is achieved when the two terms in the minimum are equal. As a result, we have 
		\begin{equation} \label{denom}
			V_A(\alpha, \beta)=p \, I_A^*+ q \, J_A^*= r(\alpha, \beta) \; q J_A^*
		\end{equation}
		with $p$ and $q$ equal to $x_A$ and $y_A$  as in 
		\eqref{LB_gap_int_1}, but  with $r$ replaced by $r(\alpha, \beta)$. 
		
		\item   Otherwise, the second term in the minimum is smaller  and the first equality in \eqref{denom} holds with $p$ and $q$ equal to $x_A$ and $y_A$  as in 
		\eqref{LB_gap_int_2}, but  with $r$ replaced by $r(\alpha, \beta)$. 
	\end{itemize}
	Therefore, letting $\alpha$ and 
	$ \beta $  go to 0  in \eqref{V_A_bou1} according to \eqref{r}, and recalling \eqref{phi}-\eqref{phi_ratio},  proves the asymptotic lower bounds in both  \eqref{LB_gap_int_1} and \eqref{LB_gap_int_2}.

	

\end{IEEEproof}

\section{} \label{app:upper_bound} 
In this Appendix we  prove Theorems \ref{theo:AO} and  \ref{theo:prob_ao}, which  provide sufficient conditions for asymptotic optimality.  In both proofs we recall that:  $c_i^*(A)>0$ for every $i$ in $A$ (resp. $A^c$) when $x_A>0$ (resp. $y_A>0$), 
$x_A \vee y_A >0$,  $x_{A}>0$ when  $ |A|>\ell$, and  
$y_{A}>0$ when  $|A|<u$. \\

\begin{IEEEproof}[Proof of Theorem \ref{theo:AO}]  We prove the theorem first when $\ell=u$, where  $\alpha$ and $\beta$ go to 0 at arbitrary rates.  By  the asymptotic  lower bound  \eqref{LB_gap}  in Theorem \ref{theo:lowerbound}  it follows that in this case  it suffices to show
	\begin{equation} \label{upper_bound_known}
		\Exp_{A}[T^{R}] \lesssim \frac{|\log (\alpha \wedge \beta)|}{x_A \, I_A^*+ y_A \, J_A^*}.
	\end{equation}
	Then, for any $\epsilon>0$ small enough and  $c>0$  we set 
	\begin{equation}\label{L}
		L_\epsilon(c):=\max_{i \in A, j \notin A} \frac{c}{ (c_i^*(A)-\epsilon) I_i+  (c_j^*(A)-\epsilon)  J_j-\epsilon},
	\end{equation}
	and observe that 
	\begin{equation}\label{expt}
		\Exp_A[T^R] \leq L_\epsilon(c)+ \sum_{n>L_\epsilon(c)} \Pro_A(T^R>n).
	\end{equation}
	For any $n \in \bN$, by  the definition of $T^R$ in \eqref{gap}  it follows that on the event $\{T^R>n\}$ there are  
	$i \in A$ and $j \notin A$ such that  $\Lambda^R_{ij}(n) <c$,  and as a result
	\begin{align*}
		\Pro_A(T^R>n) &\leq \sum_{i\in A,  j  \notin A}  \Pro_A(\Lambda^R_{ij}(n) < c). 
	\end{align*}
	Moreover, for any  $n >L_\epsilon(c)$  and $i \in A, j \notin A$ we have 
	\begin{align}  \label{c_n}
		c<n \left(  (c_i^*(A)-\epsilon) I_i+  (c_j^*(A)-\epsilon)  J_j-\epsilon \right),
	\end{align} 
	and consequently  for every $c>0$   the series in \eqref{expt} is bounded by 
	\begin{align} 
		\sum_{i\in A,  j \notin A}  \sum_{n=1}^\infty \Pro_A \left( 
		\frac{\Lambda^R_{ij}(n)}{n}<
		(c_i^*(A)-\epsilon) I_i+  (c_j^*(A)-\epsilon)  J_j-\epsilon  \right). \label{series_bound2}
	\end{align}
	By  the assumption of the theorem and an application of  Lemma \ref{new}(iii) with  $\rho_i$ equal to  $c_i^*(A)- \epsilon$ (resp. $0$) when $x_A>0$  (resp. $x_A=0$) and $\rho_j$ equal to  $c_j^*(A)- \epsilon$ (resp. $0$) when $y_A>0$ (resp. $y_A=0$),
	it follows  that the series in   \eqref{series_bound2} converges.  
	Thus, letting first $c \to \infty$ and then $\epsilon \to 0$ in \eqref{expt}  proves that   as $c \to \infty$ we have 
	\begin{equation} \label{upper_bound_known2}
		\Exp_{A}[T^{R}] \lesssim  \max_{i \in A, j \notin A}\frac{c}{c_i^*(A) I_i +  c_j^*(A) \, J_j}.
	\end{equation}
	In view of \eqref{c_star} and the selection of threshold $c$ according to \eqref{thresholds_gap}, this proves \eqref{upper_bound_known}.\\
	
	We next consider the case  $\ell<u$,   where we let  $\alpha, \beta \to 0$ so that \eqref{r} holds for some $r \in (0, \infty)$. We prove the result  when   $\ell \leq  |A| <u$, as the proof when  $\ell < |A|\leq u$ is similar.  Thus, in what follows,  $\ell \leq  |A| <u$, and as a result  $y_A>0$ and  $c_{j}^*(A) >0$ for every $j \notin A$.  By the universal  asymptotic lower bounds \eqref{LB_int} and \eqref{LB_00}   it follows that  when   either $|A|>\ell$ or   $|A|=\ell=0$,  it suffices to show that
	\begin{equation}  \label{upper_bound_unknown}
		\Exp_{A}[T^{R}] \lesssim \frac{|\log \beta|}{y_A\,  J^{*}_{A} }.
	\end{equation} 
	On the other hand, by the universal  asymptotic lower bounds \eqref{LB_gap_int_1}  and \eqref{LB_gap_int_2}    it follows that  when $|A|=\ell >0$, it suffices to show that 
	\begin{equation} \label{upper_bound_unknown2}
		\Exp_A[T^R] \lesssim  \frac{|\log \beta|}{ y_A  J_A^* }  \bigvee   \frac{ |\log \alpha|}{ x_A I_A^*+ y_A  J_A^*}.
	\end{equation} 
	When in particular  $r \leq 1$, the maximum is attained strictly by the first term, when  $r>1$,   $z_{A}<1$ and $K>\hat{K}_{A}+z_{A}\check{K}_{A}$, the maximum  is attained strictly by the second term, whereas in all other cases the two terms are equal to a first-order asymptotic approximation. 
	

	We start by  proving \eqref{upper_bound_unknown} when $\ell<|A|<u$. In this case  we also have $x_A>0$, and consequently 
	$c_{i}^*(A) >0$ for every $i \in A$. Then, for $\epsilon >0$ small   enough and $a,\, b>0$ we set 
	\begin{equation}
		N_\epsilon(a,b):= \max_{j \notin A,\, i \in A} \left\{ \frac{a}{
			(c_j^*(A)-\epsilon)  J_j-\epsilon}  ,\,   \frac{b}{
			(c_i^*(A)-\epsilon)  I_i-\epsilon}    \right\},
	\end{equation}
	and observe that 
	\begin{equation}\label{expt20}
		\Exp_A[T^R] \leq  N_\epsilon(a,b)+ \sum_{n> N_\epsilon(a,b)} \Pro_A(T^R>n).
	\end{equation}
	By the definition of $T^R$ in \eqref{gap_intersection} it follows that, for any $n \in \bN$, on the event $\{T^R>n\}$ there is either a    $j \notin A$ such that $\Lambda^R_j(n)>-a$, or an  $i \in A$ such that $\Lambda^R_i(n)<b$.  As a result, by the  union bound we obtain
	\begin{align*}
		\Pro_A(T^R>n) &\leq \sum_{j \notin A}  \Pro_A(-\Lambda^R_{j}(n) <a) + \sum_{i \in A}  \Pro_A(\Lambda^R_{i}(n) <b).
	\end{align*}
	Moreover, for any $n >  N_\epsilon(a,b)$  and  any  $i \in A$, $ j \notin A$ we have $a < n ((c_j^*(A)-\epsilon)  J_j-\epsilon)$ and $b < n ((c_i^*(A)-\epsilon)  I_i-\epsilon)$, which implies that the series in  \eqref{expt20} is bounded by     
	\begin{align}\label{series_bound10}
		\begin{split}
			&\sum_{j \notin A} \sum_{n=1}^\infty \Pro_A \left( 
			-\frac{1}{n} \Lambda^R_j(n)<
			(c_j^*(A)-\epsilon)  J_j -\epsilon  \right) \\
			&+  \sum_{i \in A} \sum_{n=1}^\infty \Pro_A \left( 
			\frac{1}{n} \Lambda^R_i(n)<
			(c_i^*(A)-\epsilon)  I_i -\epsilon  \right) .
		\end{split}
	\end{align}
	By  the assumption of the theorem and an application of  Lemma \ref{new}(i) with  $\rho_i=c_i^*(A)- \epsilon$ and of  Lemma \ref{new}(ii) with  $\rho_j=c_j^*(A)- \epsilon$ 
	it follows  that \eqref{series_bound10} converges. Thus, letting first $a,\, b \to \infty$ and then $\epsilon \to 0$ in \eqref{expt20} proves that as $a,\, b \to \infty$  
	\begin{equation*}
		\Exp_A[T^R] \lesssim  \max_{j \notin A,\, i \in A}  \left\{\frac{a}{ c_j^*(A)   J_j },\, \frac{b}{ c_i^*(A)   I_i } \right\}.
	\end{equation*}
	In view of  \eqref{c_star} and  the selection of thresholds $a, b$ according to  \eqref{thresholds_gi}, this implies that    \begin{equation} 
		\Exp_{A}[T^{R}] \lesssim \frac{|\log \beta|}{y_A\,  J^{*}_{A} } \sim \frac{|\log \alpha|}{x_A\,  I^{*}_{A} },
	\end{equation}
	and proves \eqref{upper_bound_unknown}. \\

	The proof when $|A|=\ell=0$, in which case $x_A = 0$, is similar, with the difference that we use
	\begin{equation} \label{N}
		N_\epsilon(a):= \max_{j \notin A} \left\{ \frac{a}{
			(c_j^*(A)-\epsilon)  J_j-\epsilon}   \right\}
	\end{equation}
	in the place of  $N_\epsilon(a,c)$, and apply only Lemma \ref{new}(i). \\ 

	It remains to show that \eqref{upper_bound_unknown2} holds when $|A|=\ell>0$, in which case $x_A$ is not always positive.  We recall the definitions of $L_\epsilon(c)$ and $N_\epsilon(a)$ in \eqref{L} and  \eqref{N} and observe that for any  $\epsilon >0$  small enough and $a, c>0$ we  have 
	\begin{equation}\label{expt3}
		\Exp_A[T^R] \leq L_\epsilon(c)  \vee N_\epsilon(a)+ \sum_{n>L_\epsilon(c) \vee N_\epsilon(a)} \Pro_A(T^R>n).
	\end{equation}
	By the definition of $T^R$ in \eqref{gap_intersection} it follows that, for any $n \in \bN$, on the event $\{T^R>n\}$ there are   either  $i \in A$ and $j \notin A$ such that  $\Lambda^R_{ij}(n) <c$ or  $j \notin A$ such that $\Lambda^R_j(n)>-a$, and as a result
	\begin{align*}
		\Pro_A(T^R>n) &\leq \sum_{i\in A,  j \notin A}  \Pro_A(\Lambda^R_{ij}(n) <c) +  \sum_{ j \notin A} \Pro_A(\Lambda^R_j(n) >- a ).
	\end{align*}
	Following similar steps as in the previous cases,  applying in particular  
	Lemma \ref{new}(ii)  with $\rho_j=c_j^*(A)-\epsilon$ and   Lemma \ref{new}(iii)  with $\rho_j=c_j^*(A)-\epsilon$ and 
	$\rho_i$ equal to $c_i^*(A)-\epsilon$ (resp. $0$) when $x_A>0$ (resp. $x_A=0$), 
	we conclude that as $a,c \to \infty$ 
	\begin{equation*}
		\Exp_A[T^R] \lesssim  \max_{i \in A, j \notin A} \left\{ \frac{a}{ c_j^*(A)   J_j }  \bigvee   \frac{c}{ c_i^*(A) I_i+  c_j^*(A)  J_j}\right\}.
	\end{equation*}
	In view of  \eqref{c_star} and  the selection of thresholds $a,c$ according to  \eqref{thresholds_gi},  this proves 
	\eqref{upper_bound_unknown2}. \\
\end{IEEEproof}

\begin{IEEEproof}[Proof of Theorem \ref{theo:prob_ao}]  Fix $A \in \cP_{\ell, u}$ and  a probabilistic sampling rule $R$  that satisfies \eqref{system}. Since  $c_i^*(A)>0$ for every $i$ in $A$ (resp. $A^c$) when $x_A>0$ (resp. $y_A>0$),  $x_A \vee y_A >0$, $x_{A}>0$ when  $ |A|>\ell$, and  
	$y_{A}>0$ when  $|A|<u$, the exponentially consistency of $R$  under $\Pro_A$ follows  by an application of Theorem \ref{theo:prob_cons}.  To establish its asymptotic optimality, by Theorem \ref{theo:AO} it follows that  it suffices to show that  $\Pro_A(\pi_i^R(n)<\rho)$ is  an exponentially decaying sequence for  every $\rho \in (0, c_{i}^*(A))$ and $i \in [M]$ such that $c_{i}^*(A)>0$.
	Fix such $i$ and $\rho$. Then,  there is an  $\epsilon>0$ such that 
	\begin{align} \label{rho_eps}
		\rho+\epsilon< c_{i}^*(A).
	\end{align}  
	
	By the definition of a  probabilistic rule (recall \eqref{q}), $R(n+1)$ is conditionally independent of $\mathcal{F}_{n}^{R}$ given $\Delta_{n}^{R}$ and its conditional distribution, $q^R$, does not depend on  $n$.  Thus, by \cite[Prop. 6.13]{kallenberg2002foundations} there is a measurable function 
	$h: \cP_{\ell, u} \times [0,1] \to 2^{[M]}$, which does not depend on $n$,  such that
	\begin{align*} 
		R(n+1) = h \left(\Delta_{n}^R, Z_{0}(n) \right), \quad n \in \bN,
	\end{align*}  
	where $\{Z_{0}(n), n \in \bN\}$ is a sequence of iid random variables, uniformly distributed in $(0,1)$.  Consequently,  there is a  measurable function  
	$h_i: \cP_{\ell, u} \times [0,1] \to \{0,1\}$ such that 
	\begin{align}  \label{repre_R_i}
		R_i(n+1) = h_i \left(\Delta_{n}^R, Z_0(n) \right), \quad n \in \bN.
	\end{align}

	Then, for every $n \in \bN$ we have 
	\begin{align*}
		\{ \pi_i^R(n)<\rho \} 
		&= \left\{ \sum_{m=1}^n h_i \left(A, Z_0(m-1) \right)  +    \sum_{m=1}^n (R_{i}(m)-h_i (A, Z_0(m-1) )  <n (\rho+\epsilon)- n \epsilon   \right\},
	\end{align*}
	and as a result
	\begin{align}\label{total}
		\begin{split}
			\Pro_A(\pi_i^R(n)<\rho)
			&\leq \Pro_A\left(  \sum_{m=1}^n h_i (A, Z_0(m-1))  <n (\rho+\epsilon) \right)\\
			&+\Pro_A\left( \sum_{m=1}^n ( R_{i}(m)-h_i (A, Z_0(m-1)) ) < -n \epsilon   \right). 
		\end{split}
	\end{align}
	From  \eqref{c} and \eqref{repre_R_i} it follows that 
	$\{h_i (A, Z_0(n-1)), n \in \bN\}$ is a sequence of iid Bernoulli random variables  with parameter  $c^R_{i}(A)$, whereas  by \eqref{system} and \eqref{rho_eps} it follows that $\rho+\epsilon<c^R_i(A)$. Therefore,  by the   Chernoff bound we  conclude that the first term in the upper bound in  \eqref{total}  is exponentially  decaying.   The second  term is bounded as follows
	\begin{align} \label{sec}
		&\Pro_A\left( \sum_{m=1}^n ( R_{i}(m)-h_i (A, Z_0(m-1)) ) < -n \epsilon   \right) \nonumber \\
		&\leq  \Pro_{A} \left(\sum_{m=1}^n |R_i(m) -h_i (A, Z_0(m-1)) | > n \epsilon \right) \nonumber  \\
		&\leq   \Pro_{A} \left( \sigma^R_A> n \right)+  \Pro_{A} \left( \sigma^R_A> n\epsilon \right),
	\end{align}  
	where the first inequality follows from the triangle inequality and the second by an application of the total probability rule on the event $\{ \sigma^R_A \leq n \}$, in view of the fact  that 
	$$ \sigma^R_A \leq n \quad \Rightarrow \quad 
	\sum_{m=1}^n |R_i(m) -h_i (A, Z_0(m-1)) | \leq \sigma^R_A.$$
	By  the exponentially consistency of $R$, the upper bound in \eqref{sec} is  exponentially decaying, which means that the second term in the upper bound in  \eqref{total}  is also exponentially  decaying, and this completes the proof.
\end{IEEEproof}

\bibliographystyle{IEEEtran}
\bibliography{biblio_new}

\begin{thebibliography}{10}
\providecommand{\url}[1]{#1}
\csname url@samestyle\endcsname
\providecommand{\newblock}{\relax}
\providecommand{\bibinfo}[2]{#2}
\providecommand{\BIBentrySTDinterwordspacing}{\spaceskip=0pt\relax}
\providecommand{\BIBentryALTinterwordstretchfactor}{4}
\providecommand{\BIBentryALTinterwordspacing}{\spaceskip=\fontdimen2\font plus
\BIBentryALTinterwordstretchfactor\fontdimen3\font minus
  \fontdimen4\font\relax}
\providecommand{\BIBforeignlanguage}[2]{{%
\expandafter\ifx\csname l@#1\endcsname\relax
\typeout{** WARNING: IEEEtran.bst: No hyphenation pattern has been}%
\typeout{** loaded for the language `#1'. Using the pattern for}%
\typeout{** the default language instead.}%
\else
\language=\csname l@#1\endcsname
\fi
#2}}
\providecommand{\BIBdecl}{\relax}
\BIBdecl

\bibitem{Brain_2010}
J.~Stiles and T.~Jernigan, ``The basics of brain development,''
  \emph{Neuropsychology review}, vol.~20, pp. 327--48, 11 2010.

\bibitem{Fraud_2002}
\BIBentryALTinterwordspacing
R.~J. Bolton and D.~J. Hand, ``{Statistical Fraud Detection: A Review},''
  \emph{Statistical Science}, vol.~17, no.~3, pp. 235 -- 255, 2002. [Online].
  Available: \url{https://doi.org/10.1214/ss/1042727940}
\BIBentrySTDinterwordspacing

\bibitem{De_and_Baron_2012a}
S.~K. De and M.~Baron, ``Sequential bonferroni methods for multiple hypothesis
  testing with strong control of family-wise error rates i and ii,''
  \emph{Sequential Analysis}, vol.~31, no.~2, pp. 238--262, 2012.

\bibitem{De_and_Baron_2012b}
------, ``Step-up and step-down methods for testing multiple hypotheses in
  sequential experiments,'' \emph{Journal of Statistical Planning and
  Inference}, vol. 142, no.~7, pp. 2059--2070, 2012.

\bibitem{Bartroff_and_Song_2014}
J.~Bartroff and J.~Song, ``Sequential tests of multiple hypotheses controlling
  type i and ii familywise error rates,'' \emph{Journal of statistical planning
  and inference}, vol. 153, pp. 100--114, 2014.

\bibitem{Song_and_Fellouris_2016}
Y.~Song and G.~Fellouris, ``Asymptotically optimal, sequential, multiple
  testing procedures with prior information on the number of signals,''
  \emph{Electronic Journal of Statistics}, vol.~11, 03 2016.

\bibitem{Song_and_Fellouris_2019}
\BIBentryALTinterwordspacing
------, ``{Sequential multiple testing with generalized error control: An
  asymptotic optimality theory},'' \emph{The Annals of Statistics}, vol.~47,
  no.~3, pp. 1776 -- 1803, 2019. [Online]. Available:
  \url{https://doi.org/10.1214/18-AOS1737}
\BIBentrySTDinterwordspacing

\bibitem{Bartroff_and_Song_2020}
\BIBentryALTinterwordspacing
J.~Bartroff and J.~Song, ``Sequential tests of multiple hypotheses controlling
  false discovery and nondiscovery rates,'' \emph{Sequential Analysis},
  vol.~39, no.~1, pp. 65--91, 2020. [Online]. Available:
  \url{https://doi.org/10.1080/07474946.2020.1726686}
\BIBentrySTDinterwordspacing

\bibitem{Zigangirov_1966}
\BIBentryALTinterwordspacing
K.~S. Zigangirov, ``On a problem in optimal scanning,'' \emph{Theory of
  Probability \& Its Applications}, vol.~11, no.~2, pp. 294--298, 1966.
  [Online]. Available: \url{https://doi.org/10.1137/1111025}
\BIBentrySTDinterwordspacing

\bibitem{Klimko1975OptimalSS}
E.~M. Klimko and J.~Yackel, ``Optimal search strategies for wien{\'e}r
  processes,'' \emph{Stochastic Processes and their Applications}, vol.~3, pp.
  19--33, 1975.

\bibitem{Dragalin_1996}
V.~Dragalin, ``A simple and effective scanning rule for a multi-channel
  system,'' \emph{Metrika}, vol.~43, pp. 165--182, 02 1996.

\bibitem{Cohen2015active}
K.~Cohen and Q.~Zhao, ``Asymptotically optimal anomaly detection via sequential
  testing,'' \emph{IEEE Transactions on Signal Processing}, vol.~63, no.~11,
  pp. 2929--2941, 2015.

\bibitem{huang2017active}
B.~Huang, K.~Cohen, and Q.~Zhao, ``Active anomaly detection in heterogeneous
  processes,'' \emph{IEEE Transactions on Information Theory}, vol.~65, no.~4,
  pp. 2284--2301, 2019.

\bibitem{oddball_2018}
N.~K. Vaidhiyan and R.~Sundaresan, ``Learning to detect an oddball target,''
  \emph{IEEE Transactions on Information Theory}, vol.~64, no.~2, pp. 831--852,
  2018.

\bibitem{Cohen2019nonlinearcost}
A.~Gurevich, K.~Cohen, and Q.~Zhao, ``Sequential anomaly detection under a
  nonlinear system cost,'' \emph{IEEE Transactions on Signal Processing},
  vol.~67, no.~14, pp. 3689--3703, 2019.

\bibitem{Tsopela_2019}
A.~Tsopelakos, G.~Fellouris, and V.~V. Veeravalli, ``Sequential anomaly
  detection with observation control,'' in \emph{2019 IEEE International
  Symposium on Information Theory (ISIT)}, 2019, pp. 2389--2393.

\bibitem{Cohen2020composite}
B.~Hemo, T.~Gafni, K.~Cohen, and Q.~Zhao, ``Searching for anomalies over
  composite hypotheses,'' \emph{IEEE Transactions on Signal Processing},
  vol.~68, pp. 1181--1196, 2020.

\bibitem{chernoff1959}
\BIBentryALTinterwordspacing
H.~Chernoff, ``Sequential design of experiments,'' \emph{Ann. Math. Statist.},
  vol.~30, no.~3, pp. 755--770, 09 1959. [Online]. Available:
  \url{http://dx.doi.org/10.1214/aoms/1177706205}
\BIBentrySTDinterwordspacing

\bibitem{albert1961}
\BIBentryALTinterwordspacing
A.~E. Albert, ``{The Sequential Design of Experiments for Infinitely Many
  States of Nature},'' \emph{The Annals of Mathematical Statistics}, vol.~32,
  no.~3, pp. 774 -- 799, 1961. [Online]. Available:
  \url{https://doi.org/10.1214/aoms/1177704973}
\BIBentrySTDinterwordspacing

\bibitem{Bessler1960_I}
S.~A. Bessler, ``Theory and applications of the sequential design of
  experiments, k-actions and infinitely many experiments, part i theory.''
  Department of Statistics, Stanford University, Technical Report~55, 1960.

\bibitem{Bessler1960_II}
------, ``Theory and applications of the sequential design of experiments,
  k-actions and infinitely many experiments, part ii applications.'' Department
  of Statistics, Stanford University, Technical Report~56, 1960.

\bibitem{Kiefer_Sacks_1963}
\BIBentryALTinterwordspacing
J.~Kiefer and J.~Sacks, ``{Asymptotically Optimum Sequential Inference and
  Design},'' \emph{The Annals of Mathematical Statistics}, vol.~34, no.~3, pp.
  705 -- 750, 1963. [Online]. Available:
  \url{https://doi.org/10.1214/aoms/1177704000}
\BIBentrySTDinterwordspacing

\bibitem{Keener_1984}
\BIBentryALTinterwordspacing
R.~Keener, ``{Second Order Efficiency in the Sequential Design of
  Experiments},'' \emph{The Annals of Statistics}, vol.~12, no.~2, pp. 510 --
  532, 1984. [Online]. Available: \url{https://doi.org/10.1214/aos/1176346503}
\BIBentrySTDinterwordspacing

\bibitem{Lalley_Lorden_1986}
\BIBentryALTinterwordspacing
S.~P. Lalley and G.~Lorden, ``{A Control Problem Arising in the Sequential
  Design of Experiments},'' \emph{The Annals of Probability}, vol.~14, no.~1,
  pp. 136 -- 172, 1986. [Online]. Available:
  \url{https://doi.org/10.1214/aop/1176992620}
\BIBentrySTDinterwordspacing

\bibitem{nitinawarat_controlled_2013}
S.~Nitinawarat, G.~K. Atia, and V.~V. Veeravalli, ``Controlled sensing for
  multihypothesis testing,'' \emph{IEEE Transactions on Automatic Control},
  vol.~58, no.~10, pp. 2451--2464, 2013.

\bibitem{naghshvar2013active}
M.~Naghshvar and T.~Javidi, ``Active sequential hypothesis testing,'' \emph{The
  Annals of Statistics}, vol.~41, no.~6, pp. 2703--2738, 2013.

\bibitem{nitinawarat_controlled_2015}
\BIBentryALTinterwordspacing
S.~Nitinawarat and V.~V. Veeravalli, ``Controlled sensing for sequential
  multihypothesis testing with controlled markovian observations and
  non-uniform control cost,'' \emph{Sequential Analysis}, vol.~34, no.~1, pp.
  1--24, 2015. [Online]. Available:
  \url{https://doi.org/10.1080/07474946.2014.961864}
\BIBentrySTDinterwordspacing

\bibitem{Aditya_2021}
\BIBentryALTinterwordspacing
A.~Deshmukh, V.~V. Veeravalli, and S.~Bhashyam, ``Sequential controlled sensing
  for composite multihypothesis testing,'' \emph{Sequential Analysis}, vol.~40,
  no.~2, pp. 259--289, 2021. [Online]. Available:
  \url{https://doi.org/10.1080/07474946.2021.1912525}
\BIBentrySTDinterwordspacing

\bibitem{garivier2016optimal}
A.~Garivier and E.~Kaufmann, ``Optimal best arm identification with fixed
  confidence,'' in \emph{Conference on Learning Theory}.\hskip 1em plus 0.5em
  minus 0.4em\relax PMLR, 2016, pp. 998--1027.

\bibitem{resnick1992adventures}
S.~I. Resnick, \emph{Adventures in stochastic processes}.\hskip 1em plus 0.5em
  minus 0.4em\relax Springer Science \& Business Media, 1992.

\bibitem{hommel1988controlled}
G.~Hommel and T.~Hoffmann, ``Controlled uncertainty,'' in \emph{Multiple
  Hypothesenpr{\"u}fung/Multiple Hypotheses Testing}.\hskip 1em plus 0.5em
  minus 0.4em\relax Springer, 1988, pp. 154--161.

\bibitem{LehmannRomano2005}
\BIBentryALTinterwordspacing
E.~L. Lehmann and J.~P. Romano, ``Generalizations of the familywise error
  rate,'' \emph{Ann. Statist.}, vol.~33, no.~3, pp. 1138--1154, 06 2005.
  [Online]. Available: \url{http://dx.doi.org/10.1214/009053605000000084}
\BIBentrySTDinterwordspacing

\bibitem{Tsopela_2020}
A.~Tsopelakos and G.~Fellouris, ``Sequential anomaly detection with observation
  control under a generalized error metric,'' in \emph{2020 IEEE International
  Symposium on Information Theory (ISIT)}, 2020, pp. 1165--1170.

\bibitem{heydari2016quickest}
J.~Heydari, A.~Tajer, and H.~V. Poor, ``Quickest linear search over correlated
  sequences,'' \emph{IEEE Transactions on Information Theory}, vol.~62, no.~10,
  pp. 5786--5808, 2016.

\bibitem{rag_sas13}
\BIBentryALTinterwordspacing
M.~Raginsky and I.~Sason, ``Concentration of measure inequalities in
  information theory, communications, and coding,'' \emph{Foundations and
  Trends® in Communications and Information Theory}, vol.~10, no. 1-2, pp.
  1--246, 2013. [Online]. Available: \url{http://dx.doi.org/10.1561/0100000064}
\BIBentrySTDinterwordspacing

\bibitem{wald1945sequential}
A.~Wald, ``Sequential tests of statistical hypotheses,'' \emph{The Annals of
  Mathematical Statistics}, vol.~16, no.~2, pp. 117--186, 1945.

\bibitem{tartakovsky_book_2014}
A.~Tartakovsky, I.~Nikiforov, and M.~Basseville, \emph{Sequential analysis:
  Hypothesis testing and changepoint detection}.\hskip 1em plus 0.5em minus
  0.4em\relax CRC press, 2015.

\bibitem{Chow_Teicher2012}
\BIBentryALTinterwordspacing
Y.~Chow and H.~Teicher, \emph{Probability Theory: Independence,
  Interchangeability, Martingales}, ser. Springer Texts in Statistics.\hskip
  1em plus 0.5em minus 0.4em\relax Springer New York, 2012. [Online].
  Available: \url{https://books.google.com/books?id=213dBwAAQBAJ}
\BIBentrySTDinterwordspacing

\bibitem{kallenberg2002foundations}
\BIBentryALTinterwordspacing
O.~Kallenberg, \emph{Foundations of Modern Probability}, ser. Probability and
  Its Applications.\hskip 1em plus 0.5em minus 0.4em\relax Springer New York,
  2002. [Online]. Available:
  \url{https://books.google.com/books?id=TBgFslMy8V4C}
\BIBentrySTDinterwordspacing

\end{thebibliography}

\end{document}